\documentclass[11pt]{amsart}

\usepackage[all]{xy}
\usepackage{amscd}
\usepackage{mathrsfs}
\usepackage{amssymb}
\usepackage{amsthm}
\usepackage{amsmath}

\newtheorem{thm}{Theorem}[section]
\newtheorem{prop}[thm]{Proposition}
\newtheorem{lem}[thm]{Lemma}
\newtheorem{cor}[thm]{Corollary}
\newtheorem{conj}[thm]{Conjecture}

\newtheorem{ithm}{Theorem}

\newtheorem{icor}[ithm]{Corollary}

\newcommand{\onto}{\rightarrow \hspace{-.86em} \rightarrow}

\newcommand{\ra}{\rightarrow}
\newcommand{\lra}{\longrightarrow}

\newcommand{\into}{\hookrightarrow}

\newcommand{\iso}{\stackrel{\sim}{\ra}}
\newcommand{\liso}{\stackrel{\sim}{\lra}}

\newcommand{\pfbegin}{{{\em Proof:}\;}}
\newcommand{\pfend}{$\Box$ \medskip}

\newlength{\ownl}

\newcommand{\norm}{{\mbox{\bf N}}}
\newcommand{\ndiv}{{\mbox{$\not| $}}\ }

\newcommand{\Ma}{{\operatorname{Ma}}}

\newcommand{\Art}{{\operatorname{Art}}}
\DeclareMathOperator{\Aut}{\mathrm{Aut}}
\DeclareMathOperator{\BC}{\mathrm{BC}}

\newcommand{\End}{{\operatorname{End}\,}}

\newcommand{\Fil}{{\operatorname{Fil}\,}}

\newcommand{\Frob}{{\operatorname{Frob}}}
\newcommand{\Gal}{{\operatorname{Gal}\,}}
\newcommand{\Groth}{{\operatorname{Groth}\,}}
\newcommand{\Hom}{{\operatorname{Hom}\,}}

\newcommand{\Iw}{{\operatorname{Iw}}}

\newcommand{\Lie}{{\operatorname{Lie}\,}}

\newcommand{\Red}{{\operatorname{Red}}}

\newcommand{\Spec}{{\operatorname{Spec}\,}}
\newcommand{\Spf}{{\operatorname{Spf}\,}}
\newcommand{\Spp}{{\operatorname{Sp}}}

\newcommand{\HT}{{\operatorname{HT}}}
\newcommand{\semit}{{\operatorname{st}}}

\newcommand{\gr}{{\operatorname{gr}\,}}
\DeclareMathOperator{\pr}{\mathrm{pr}}
\newcommand{\rec}{{\operatorname{rec}}}
\DeclareMathOperator{\tr}{\mathrm{tr}}

\newcommand{\vol}{{\operatorname{vol}\,}}

\newcommand{\nind}{{\operatorname{n-Ind}\,}}

\newcommand{\ab}{{\operatorname{ab}}}

\newcommand{\dr}{{\operatorname{DR}}}
\newcommand{\et}{{\operatorname{et}}}

\newcommand{\WD}{{\operatorname{WD}}}

\newcommand{\op}{{\operatorname{op}}}

\newcommand{\semis}{{\operatorname{ss}}}
\newcommand{\Fsemis}{{F{\operatorname{-ss}}}}

\newcommand{\A}{{\Bbb{A}}}

\newcommand{\C}{{\Bbb{C}}}

\newcommand{\F}{{\Bbb{F}}}
\newcommand{\G}{{\Bbb{G}}}

\newcommand{\Q}{{\Bbb{Q}}}
\newcommand{\R}{{\Bbb{R}}}

\newcommand{\Z}{{\Bbb{Z}}}

\newcommand{\CA}{{\mathcal{A}}}

\newcommand{\CC}{{\mathcal{C}}}

\newcommand{\CF}{{\mathcal{F}}}
\newcommand{\CG}{{\mathcal{G}}}
\newcommand{\CH}{{\mathcal{H}}}

\newcommand{\CK}{{\mathcal{K}}}
\newcommand{\CL}{{\mathcal{L}}}

\newcommand{\CO}{{\mathcal{O}}}

\newcommand{\gX}{{\frak{X}}}

\newcommand{\gm}{{\frak{m}}}

\newcommand{\bareta}{{\overline{\eta}}}

\newcommand{\ol}{\overline}

\def\invlim#1{\lim\limits_{\substack{\longleftarrow\\#1}}}
\def\dirlim#1{\lim\limits_{\substack{\longrightarrow\\#1}}}

\title{Compatibility of local and global Langlands correspondences}

\author{Richard Taylor}\thanks{The first author was partially supported by NSF grant number DMS-0100090}
\author{Teruyoshi Yoshida}
\address{Harvard University, Department of Mathematics, 1 Oxford Street, Cambridge, MA 02138, USA}
\email{rtaylor\char`\@math.harvard.edu}
\email{yoshida\char`\@math.harvard.edu}

\date{\today}

\begin{document}
\baselineskip 14pt

\maketitle

\begin{abstract}
We prove the compatibility of local and global Langlands correspondences for
$GL_n$, which was proved up to semisimplification in \cite{ht}. More precisely,
for the $n$-dimensional $l$-adic representation $R_l(\Pi)$ of the Galois group
of an imaginary CM-field $L$ attached to a conjugate self-dual regular algebraic cuspidal
automorphic representation $\Pi$ of $GL_n(\A_L)$, which is square integrable at some finite
place, we show that Frobenius semisimplification of the restriction of
$R_l(\Pi)$ to the decomposition group of a place $v$ of $L$ not dividing $l$
corresponds to $\Pi_v$ by the local Langlands correspondence. If $\Pi_v$ is
square integrable for some finite place $v \ndiv l$ we deduce that
$R_l(\Pi)$
is irreducible. We also obtain conditional results in the case $v|l$.
\end{abstract}

\section*{Introduction}

This paper is a continuation of \cite{ht}. Let $L$ be an imaginary CM field
and let $\Pi$ be a regular algebraic cuspidal automorphic representation of
$GL_n(\A_L)$ which is conjugate self-dual ($\Pi \circ c \cong \Pi^\vee$) and
square integrable at some finite place. In \cite{ht} it is explained how to
attach to $\Pi$ and an arbitrary rational prime $l$ (and an isomorphism $\imath:\Q_l^{ac} \iso \C$) a continuous semisimple representation
\[ R_l(\Pi): \Gal(L^{ac}/L) \lra GL_n(\Q_l^{ac}) \]
which is characterised as follows. For every finite place $v$ of $L$ not
dividing $l$ 
\[ \imath R_l(\Pi)|_{W_{L_v}}^\semis = \rec\bigl(\Pi_v^\vee
|\det|^{\frac{1-n}{2}}\bigr)^\semis, \]
where $\rec$ denotes the local Langlands correspondence and $\semis$ denotes
the semisimplification (see \cite{ht} for details). In \cite{ht} it is also
shown that $\Pi_v$ is tempered for all finite places $v$.

In this paper we strengthen this result to completely identify
$R_l(\Pi)|_{I_v}$ for $v \ndiv l$. In particular, we prove the following
theorem.
\begin{ithm}
If $v \ndiv l$ then the Frobenius
semisimplification of $R_l(\Pi)|_{W_{L_v}}$ is the $l$-adic representation
attached to $\imath^{-1} \rec\bigl(\Pi_v^\vee |\det|^{\frac{1-n}{2}}\bigr)$.
\end{ithm}
As $R_l(\Pi)$ is semisimple and $\rec\bigl(\Pi_v^\vee
|\det|^{\frac{1-n}{2}}\bigr)$ is indecomposable if $\Pi_v$ is square
integrable, we obtain the following corollary.
\begin{icor}
If $\Pi_v$ is square integrable at a finite place $v\ndiv l$, then the representation $R_l(\Pi)$ is irreducible.
\end{icor}
We also obtain some results in the case $v|l$ which we will describe in section one.

Using base change it is easy to reduce to the case that $\Pi_v$ has an Iwahori
fixed vector. We descend $\Pi$ to an automorphic representation $\pi$ of a unitary group 
$G$ which locally at $v$ looks like $GL_n$ and at infinity looks like $U(n-1,1)
\times U(n,0)^{[L:\Q]/2-1}$. Then we realise $R_l(\Pi)$ in the cohomology of a Shimura
variety $X$ associated to $G$ with Iwahori level structure at $v$. More
precisely, for some $l$-adic sheaf $\CL$, the $\pi^p$-isotypic component of
$H^{n-1}(X, \CL)$ is, up to semisimplification and some twist, $R_l(\Pi)^a$ (for some $a \in
\Z_{>0}$). We show that $X$ has strictly semistable reduction and use the results of
\cite{ht} to calculate the cohomology of the (smooth, projective) strata of
the reduction of $X$ above $p$ as a virtual $G(\A^{\infty,p}) \times
\Frob_v^\Z$-module (where $\Frob_v$ denotes Frobenius). This description and the
temperedness of $\Pi_v$ shows that the $\pi^p$-isotypic component of the
cohomology of any strata is concentrated in the middle degree. This implies
that the $\pi^p$-isotypic component of the Rapoport-Zink weight spectral
sequence degenerates at $E_1$, which allows us to calculate the action of
inertia at $v$ on $H^{n-1}(X, \CL)$.

In the special case that $\Pi_v$ is a twist of a Steinberg representation
and $\Pi_\infty$ has trivial infinitesimal character, the above theorem
presumably follows from the results of Ito \cite{i}.

After we had posted the first version of this paper, Boyer \cite{bo} has announced an alternative proof with presumably stronger results.

{\bf Acknowledgements}. The authors are grateful to Tetsushi Ito, who
asked a very helpful question.

\section{The main theorem}

We write $F^{ac}$ for an algebraic closure of a field $F$. Let $l$ be a rational prime
and $\imath:\Q_l^{ac} \iso \C$ an isomorphism.

Suppose that $p$ is another rational prime.
Let $K/\Q_p$ be a finite extension. We will let $\CO_K$ denote the ring of
integers of $K$, $\wp_K$ the unique maximal ideal of $\CO_K$, $v_K$ the
canonical valuation $K^\times \onto \Z$, $k(v_K)$ the residue field
$\CO_K/\wp_K$ and $| \,\,\,|_K$ the absolute value normalised by $|x|_K =
(\# k(v_K))^{-v_K(x)}$. We will let $\Frob_{v_K}$ denote the geometric
Frobenius element of $\Gal(k(v_K)^{ac}/k(v_K))$. We will let $I_{v_K}$ denote
the kernel of the natural surjection $\Gal(K^{ac}/K) \onto
\Gal(k(v_K)^{ac}/k(v_K))$. We will let $W_K$ denote the preimage under
$\Gal(K^{ac}/K) \onto \Gal(k(v_K)^{ac}/k(v_K))$ of $\Frob_{v_K}^\Z$ endowed
with a topology by decreeing that $I_K$ with its usual topology is an open
subgroup of $W_K$. Local class field theory provides a canonical isomorphism
$\Art_K: K^\times \iso W_K^{ab}$, which takes uniformisers to lifts of
$\Frob_{v_K}$.

Let $\Omega$ be an algebraically closed field of characteristic $0$ and of the
same cardinality as $\C$. (Thus in fact $\Omega \cong \C$.) By a
{\em Weil-Deligne representation} of $W_K$ over $\Omega$ we mean a finite dimensional
$\Omega$-vector space $V$ together with a homomorphism $r:W_K \ra GL(V)$ with
open kernel and an element $N \in \End(V)$ which satisfies
\[ r(\sigma) N r(\sigma)^{-1} = |\Art^{-1}_K (\sigma)|_K N .\]
We sometimes denote a Weil-Deligne representation by $(V,r,N)$ or simply
$(r,N)$. For a finite extension $K'/K$ of $p$-adic fields, we define
\[ (V,r,N)|_{W_{K'}} = (V,r|_{W_{K'}}, N). \]

We call $(V,r,N)$ {\em Frobenius semisimple} if $r$ is semisimple. If $(V,r,N)$ is
any Weil-Deligne representation we define its {\em Frobenius semisimplification}
$(V,r,N)^{\Fsemis}=(V,r^{\semis},N)$ as follows. Choose a lift $\phi$ of
$\Frob_{v_K}$ to $W_K$. Let $r(\phi)=su=us$ where $s \in GL(V)$ is semisimple
and $u \in GL(V)$ is unipotent. For $n \in \Z$ and $\sigma \in I_K$ set
$r^{\semis}(\phi^n\sigma)=s^n r(\sigma)$. This is independent of the choices, 
and gives a Frobenius semisimple Weil-Deligne representation. We will also
set $(V,r,N)^\semis=(V,r^\semis,0)$.

One of the main results of \cite{ht} is that, given a choice of
$(\# k(v_K))^{1/2} \in \Omega$, there is a bijection $\rec$ (the local Langlands correspondence) from isomorphism
classes of irreducible smooth representations of $GL_n(K)$ over $\Omega$ to
isomorphism classes of $n$-dimensional Frobenius semisimple Weil-Deligne
representations of $W_K$, and that this bijection is natural in a number of
respects. (See \cite{ht} for details.)

Suppose that $l \neq p$. We will call a Weil-Deligne representation of $W_K$ over $\Q_l^{ac}$
{\em bounded} if for some (and hence all) $\sigma \in W_K-I_K$,
all the eigenvalues of $r(\sigma)$ are $l$-adic units. There is an equivalence
of categories between bounded Weil-Deligne representations of $W_K$ over
$\Q_l^{ac}$ and continuous representations of $\Gal(K^{ac}/K)$ on finite
dimensional $\Q_l^{ac}$-vector spaces as follows. Fix a lift $\phi\in W_K$ of
$\Frob_{v_K}$ and a continuous homomorphism $t:I_K \onto \Z_l$. Send a
Weil-Deligne representation $(V,r,N)$ to $(V,\rho)$, where $\rho$ is the unique
continuous representation of $\Gal(K^{ac}/K)$ on $V$ such that
\[ \rho(\phi^n \sigma) = r(\phi^n \sigma) \exp \bigl( t(\sigma) N \bigr) \]
for all $n \in \Z$ and $\sigma \in I_K$. Up to natural isomorphism this
functor is independent of the choices of $t$ and $\phi$. We will write 
$\WD(V,\rho)$ for the Weil-Deligne representation corresponding to a continuous
representation $(V,\rho)$. If $\WD(V,\rho)=(V,r,N)$, then have
$\rho|_{W_K}^\semis \cong r^{\semis}$. (See \cite{tate}, \S 4 and \cite{d},
\S 8 for details.) 

Now suppose that $l=p$. Let $K_0$ denote the maximal subfield of $K$ which is
unramified over $\Q_p$. Recall the filtered $K$-algebra $B_\dr$ and the $K_0$-algebra
$B_\semit$ with endomorphisms $N$ and $\phi$, where $\phi$ is $\Frob_{\Q_p}^{-1}$-semilinear
and $\phi N = p N \phi$. (See \cite{pp}.) A continuous representation $(V,\rho)$ of
$\Gal(K^{ac}/K)$ on a finite dimensional $\Q_l^{ac}$ vector space is called {\em de Rham}
if 
\[ (V \otimes_{\Q_l} B_{\dr})^{\Gal(K^{ac}/K)} \]
is a free over $\Q_l^{ac} \otimes_{\Q_l} K$ of rank $\dim_{\Q_l^{ac}} V$. The category
of de Rham representations is closed under tensor operations. One can also define
various invariants of de Rham representations $(V,\rho)$. Firstly, for each embedding $\tau:K \into
\Q_l^{ac}$, we will let $\HT_\tau(V,\rho)$ denote the $\dim_{\Q_l^{ac}} V$ element
multiset of integers which contains $j$ with multiplicity 
\[ \dim_{\Q_l^{ac}}  \gr^j (V \otimes_{\tau, K} B_{\dr})^{\Gal(K^{ac}/K)} = \dim_{\Q_l^{ac}}
\gr^j (V \otimes_{\Q_l} B_{\dr})^{\Gal(K^{ac}/K)} \otimes_{(\Q_l^{ac} \otimes_{\Q_l} K), 1 \otimes \tau}
\Q_l^{ac}. \]
(These are referred to as the Hodge-Tate numbers of $V$ with respect to $\tau$.)
Moreover if $(V,\rho)$ is de Rham then we can find a finite Galois extension $L/K$ such that
\[ (V \otimes_{\Q_l} B_\semit)^{\Gal(K^{ac}/L)} \]
is a free $\Q_l^{ac} \otimes_{\Q_l} L_0$ module of rank
 $\dim_{\Q_l^{ac}} V$, where $L_0/\Q_p$ is the maximal unramified
 subextension of $L/\Q_p$  (see \cite{berg}). Choose such an
extension $L$ and choose $\tau:L_0 \into \Q_l^{ac}$. Then set
\[ W=(V \otimes_{\tau,L} B_\semit)^{\Gal(K^{ac}/L)} = (V \otimes_{\Q_l} B_\semit)^{\Gal(K^{ac}/L)}
\otimes_{(\Q_l^{ac} \otimes_{\Q_l} L_0), 1 \otimes \tau} \Q_l^{ac}. \]
For $\sigma \in W_K$ lying above $\Frob_K^m$ set
\[ r(\sigma)=\rho(\sigma) \otimes (\sigma \phi^{m[k(v_K):\F_l]}) \in GL(W). \]
Finally set
\[ \WD(V,\rho)=(W,r,1 \otimes N). \]
This is a Weil-Deligne representation of $W_K$, which up to isomorphism is independent of
the choices of $L$ and $\tau$. In fact $\WD$ is a functor from de Rham representations of
$\Gal(K^{ac}/K)$ over $\Q_l^{ac}$ to Weil-Deligne representations of $W_K$ over $\Q_l^{ac}$.
(These constructions are due to Fontaine. See \cite{pp} for details.) 

In both the cases $l \neq p$ and $l=p$ the functor $\WD$ commutes with restriction to
an open subgroups and tensor operations. Recall the following standard conjecture.
\begin{conj}\label{il} Suppose that $X/K$ is a proper smooth variety purely of dimension $n$.  Suppose
also that $\sigma \in W_K$ and that $\Gamma \in CH^n(X \times_K X)$ is an algebraic
correspondence. Then the alternating sum of the trace
\[ \sum_{i=0}^{2n} (-1)^i \tr \bigl(\sigma \Gamma^*\big| \WD\bigl(H^i(X \times_K K^{ac}, \Q_l^{ac})\bigr)\bigr) \]
lies in $\Q$ and is independent of $l$. Here $\sigma$ induces the automorphism by the right action on $X\times_KK^{ac}$, and $\Gamma^*$ is the endomorphism defined by $\pr_{1,*}\circ ([\Gamma]\cup)\circ \pr_2^*$.
\end{conj}
For $l \neq p$ this conjecture was proven by T.\ Saito \cite{saito}. In the introduction to that paper he
expresses the conviction that the case $l=p$ can be handled by the same methods. This was carried
out by Ochiai \cite{och} in the case $\Gamma$ is the diagonal $X \subset X \times_K X$.

Now we turn to some global considerations. Let $L$ be a finite, imaginary CM  extension of $\Q$.
Let $c \in \Aut(L)$ denote
complex conjugation. Suppose that $\Pi$ is a cuspidal automorphic
representation of $GL_n(\A_L)$ such that
\begin{itemize}
\item $\Pi \circ c \cong \Pi^\vee$;
\item $\Pi_\infty$ has the same infinitesimal character as some algebraic
representation over $\C$ of the restriction of scalars from $L$ to $\Q$ of
$GL_n$;
\item and for some finite place $x$ of $L$ the representation $\Pi_x$ is
square integrable.
\end{itemize}
(In this paper `square integrable' (resp.\ `tempered') will mean the twist
by a character of a pre-unitary representation which is square integrable
(resp.\  tempered).)
In \cite{ht} (see theorem C in the introduction of \cite{ht}) it is shown that
there is a unique continuous semisimple representation
\[ R_{l,\imath}(\Pi)=R_l(\Pi):\Gal(L^{ac}/L) \lra GL_n(\Q_l^{ac}) \]
such that for each finite place $v\ndiv l$ of $L$ 
\[ \rec\bigl(\Pi_v^\vee
|\det|^{\frac{1-n}{2}}\bigr)^\semis=(\imath R_l(\Pi)|_{W_{L_v}}^\semis, 0). \]
Moreover it is shown that $\Pi_v$ is tempered for all finite
places $v$ of $L$, which completely determines the $N$ in $\rec\bigl(\Pi_v^\vee
|\det|^{\frac{1-n}{2}}\bigr)$ (see lemma \ref{pure}
below). If $n=1$ both these assertions are true without the
assumptions that $\Pi \circ c \cong \Pi^\vee$ and (for the first assertion) $v \ndiv l$.

The main theorem of this paper identifies
$\WD(R_l(\Pi)|_{\Gal(L_v^{ac}/L_v)})^\Fsemis$. More precisely we
prove the following.

\begin{thm}\label{mt} Keep the above notation and assumptions. Then for each
finite place $v \ndiv l$ of $L$ there is an isomorphism 
\[ \imath \WD(R_l(\Pi)|_{\Gal(L_v^{ac}/L_v)})^\Fsemis \cong \rec\bigl(\Pi_v^\vee
|\det|^{\frac{1-n}{2}}\bigr)\] 
of Weil-Deligne representations over $\C$. If conjecture \ref{il} is true then this holds
even for $v|l$. \end{thm}

As $R_l(\Pi)$ is semisimple and $\rec\bigl(\Pi_v^\vee |\det|^{\frac{1-n}{2}}
\bigr)$ is indecomposable if $\Pi_v$ is square integrable, we have the
following corollary.

\begin{cor}
If $\Pi_v$ is square integrable for a finite place $v\ndiv l$ or if conjecture \ref{il} is true, then
the representation $R_l(\Pi)$ is irreducible. (Recall that we are assuming that $\Pi_v$ is square integrable
for some finite place $v$.) 
\end{cor}

In the rest of this section we consider some generalities on Galois
representations and Weil-Deligne representations. First
consider Weil-Deligne representations over an algebraically closed field
$\Omega$ of characteristic zero and the same cardinality as $\C$.
If $(W,r)$ is a finite dimensional representation of $W_K$ with open kernel and
if $s \in \Z_{\geq 1}$ we will write $\Spp_s(W)$ for the Weil-Deligne
representation
\[ \bigl(W^s,\ r |\Art_K^{-1}|_K^{s-1} \oplus \cdots \oplus r |\Art_K^{-1}|_K \oplus r,
\ N \bigr) \]
where $N: r |\Art_K^{-1}|_K^{i-1} \iso r |\Art_K^{-1}|_K^i$ for $i=1,...,s-1$.
This defines $\Spp_s(W)$ uniquely (up to isomorphism). If $W$ is irreducible
then $\Spp_s(W)$ is indecomposable and every indecomposable Weil-Deligne
representation is of the form $\Spp_s(W)$ for a unique $s$ and a unique
irreducible $W$. If $\pi$ is an irreducible cuspidal representation of
$GL_g(K)$ then $\rec(\pi)=(r,0)$ with $r$ irreducible. Moreover for any
$s \in \Z_{\geq 1}$ we have (in the notation of section I.3 of \cite{ht})
$\rec(\Spp_s(\pi))=\Spp_s(r)$.

If $q \in \R_{>0}$, then by a {\em Weil $q$-number} we mean $\alpha \in \Q^{ac}$
such that for all $\sigma: \Q^{ac} \into \C$ we have $(\sigma \alpha)(c \sigma
\alpha)=q$. We will call a Weil-Deligne representation $(V,r,N)$ of $W_K$
{\em strictly pure of weight $k\in \R$} if for some (and hence every) lift
$\phi$ of $\Frob_{v_K}$, every eigenvalue $\alpha$ of $r(\phi)$ is a Weil
$(\# k(v_K))^k$-number. In this case we must have $N=0$. We will call $(V,r,N)$
{\em mixed} if it has an increasing filtration $\Fil^W_i$ with $\Fil^W_i V=V$
for $i \gg 0$ and $\Fil^W_i V=(0)$ for $i\ll 0$, such that the
$i$-th graded piece is strictly pure of weight $i$. If $(V,r,N)$ is mixed
then there is a unique choice of filtration $\Fil^W_i$, and
$N(\Fil^W_iV)\subset \Fil^W_{i-2}V$. Finally we will call
$(V,r,N)$ {\em pure of weight $k$} if it is mixed with all weights in $k+\Z$
and if for all $i \in \Z_{>0}$
\[ N^i: \gr^W_{k+i} V \liso \gr_{k-i}^W V. \]
If $W$ is strictly pure of weight $k$, then $\Spp_s(W)$ is pure of weight $k-(s-1)$.
(It is generally conjectured that if $X$ is a proper smooth variety over a
$p$-adic field $K$, then $\WD\bigl(H^i(X\times_K K^{ac}, \Q_l^{ac})\bigr)$ is
pure of weight $i$ in the above sense.)

\begin{lem} \label{pure} \begin{enumerate}
\item $(V,r,N)$ is pure if and only if $(V,r,N)^{\Fsemis}$ is.
\item If $L/K$ is a finite extension, then $(V,r,N)$ is pure if and only if
$(V,r,N)|_{W_L}$ is pure.
\item An irreducible smooth representation $\pi$ of
$GL_n(K)$ has $\sigma \pi$ tempered for all $\sigma:\Omega \into \C$ if and
only if $\rec(\pi)$ is pure of some weight.
\item Given $(V,r)$ with $r$ semisimple, there is, up to equivalence, at most
one choice of $N$ which makes $(V,r,N)$ pure. 
\item If $(V,r,N)$ is a Frobenius semisimple Weil-Deligne representation
which is pure of weight $k$ and if $W \subset V$ is a Weil-Deligne
subrepresentation, then the following are equivalent:
\begin{enumerate}
\item $\bigwedge^{\dim W} W$ is pure of weight $k \dim W$,
\item $W$ is pure of weight $k$,
\item $W$ is a direct summand of $V$.
\end{enumerate}
\item Suppose that $(V,r,N)$ is a Frobenius semisimple Weil-Deligne
representation which is pure of weight $k$. Suppose also that $\Fil^j V$ is a
decreasing filtration of $V$ by Weil-Deligne subrepresentations such that
$\Fil^j V=(0)$ for $j \gg 0$ and $\Fil^jV=V$ for $j\ll 0$. If for each $j$
\[ \textstyle{\bigwedge^{\dim \gr^j V} \gr^jV } \]
is pure of weight $k \dim \gr^jV$, then 
\[ V \cong \bigoplus_j \gr^j V \]
and each $\gr^j V$ is pure of weight $k$.
\end{enumerate}\end{lem}

\pfbegin
The first two parts are straightforward (using the fact that the filtration 
$\Fil^W_i$ is unique). For the third part recall that an irreducible
smooth representation $\Spp_{s_1}(\pi_1) \boxplus \cdots \boxplus
\Spp_{s_t}(\pi_t)$ (see section I.3 of \cite{ht}) is tempered if and only if
the absolute values of the central characters of the $\Spp_{s_i}(\pi_i)$ are
all equal. 

Suppose that $(V,r,N)$ is Frobenius semisimple and pure of weight $k$. As a $W_K$-module we can write
uniquely as $V = \oplus _{i \in \Z} V_i$ where $(V_i,r|_{V_i},0)$ is strictly pure of weight
$k+i$. For $i \in \Z_{\geq 0}$ let $V(i)$ denote the kernel of $N^{i+1}:V_i \ra V_{-i-2}$.
Then $N:V_{i+2} \into V_i$ is injective and $V_i=NV_{i+2} \oplus  V(i)$. Thus
\[ V = \bigoplus_{i \in \Z} \bigoplus_{j=0}^i N^jV(i), \]
and for $0 \leq j \leq i$ the map $N^j:V(i) \ra V_{i-2j}$ is injective. 
Also note that as a virtual $W_K$-module $[V(i)]=[V_i]-[V_{i+2}\otimes |\Art_K^{-1}|_K]$. Thus if $r$ is
semisimple then $(V,r)$ determines $(V,r,N)$ up to isomorphism. This establishes
the fourth part.

Now consider the fifth part. If $W$ is a direct summand it is certainly pure
of the same weight $k$ and $\bigwedge^{\dim W}W$ is then pure of weight $k\dim W$.
Conversely if $W$ is pure of weight $k$ then
\[ W = \bigoplus_{i \in \Z} \bigoplus_{j=0}^i N^jW(i), \]
where $W(i)=W \cap V(i)$. As a $W_K$-module we can decompose $V(i)=W(i) \oplus
U(i)$. Setting
\[ U = \bigoplus_{i \in \Z} \bigoplus_{j=0}^i N^jU(i), \]
we see that $V=W \oplus U$ as Weil-Deligne representations.
Now suppose only that $\bigwedge^{\dim W}W$ is pure of weight $k \dim W$. Write
\[ W \cong \bigoplus_j \Spp_{s_j}(X_j) \]
where each $X_j$ is strictly pure of some weight $k+k_j+(s_j-1)$. Then,
looking at highest exterior powers, we see that $\sum_j k_j(\dim \Spp_{s_j}(X_j))=0$. On the other
hand as $V$ is pure we see that $k_j \leq 0$ for all $j$. We conclude that
$k_j=0$ for all $j$ and hence that $W$ is pure of weight $k$. 

The final part follows from the fifth part by a simple inductive argument.
\pfend

In view of parts (3) and (4) this lemma, theorem \ref{mt} will follow from the following two results,
which we prove in the rest of this paper. (Recall that for each $\sigma \in \Aut(\C)$ the representation
$\sigma \Pi$ is again cuspidal automorphic \cite{cloz} and hence for each finite place $v$ the
representation $\sigma \Pi_v$ is tempered.)

\begin{thm}\label{mt2} Keep the notation and assumptions of theorem \ref{mt}. Then for each
finite place $v$ of $L$ 
\[ \WD(R_l(\Pi)|_{\Gal(L_v^{ac}/L_v)}) \]
is pure.
\end{thm}

\begin{prop}\label{il2} Keep the notation and assumptions of theorem \ref{mt}. Assume also that
conjecture \ref{il} is true. Suppose finally that $l'$ is a rational prime and that $\imath':\Q_{l'}^{ac} \iso \C$.
Then for each finite place $v$ of $L$ 
\[ \imath \WD(R_{l,\imath}(\Pi)|_{\Gal(L_v^{ac}/L_v)})^\semis \cong  \imath' \WD(R_{l',\imath'}(\Pi)|_{\Gal(L_v^{ac}/
L_v)})^\semis. \]
\end{prop}

Now let $L$ denote a number field. Write $|\,\,\,|_L$ for
\[ \prod_x |\,\,\,|_{L_x}: \A_L^\times/L^\times \lra \R_{>0}^\times, \]
and write $\Art_L$ for
\[ \prod_x \Art_{L_x}: \A_L^\times /L^\times \onto \Gal(L^{ac}/L)^\ab. \]
We will call a continuous representation 
\[ R: \Gal(L^{ac}/L) \lra GL_n(\Q_l^{ac}) \]
{\em pure of weight $k$} if for all but finitely many finite places $x$ of
$L$ the representation $R$ is unramified at $x$ and every eigenvalue $\alpha$
of $R(\Frob_x)$ is a Weil $(\# k(x))^k$-number.

If
\[ R: \Gal(L^{ac}/L) \lra GL_1(\Q_l^{ac}) \]
is de Rham at all places $v|l$ of $L$, then there exist
\begin{itemize}
\item a CM (possibly totally real) field $L_0 \subset L$;
\item an integer $k$;
\item integers $n(\tau)$ for each $\tau\colon L_0 \into \Q_l^{ac}$, such that $n(\tau)+n(\tau c)=k$
for all $\tau$;
\item and a continuous character $\chi: \A_L^\times \ra (\Q^{ac})^\times \subset (\Q_l^{ac})^\times$ such that
\[ \chi|_{L^\times} = \prod_{\tau:L \into \Q_l^{ac}} \tau^{n(\tau|_{L_0})}; \]
\end{itemize}
with the following properties:
\begin{itemize}
\item $\displaystyle R \circ \Art_L (x) = \chi(x) \prod_{\tau:L \into \Q_l^{ac}} \tau(x_l)^{-n(\tau|_{L_0})}$;
\item and for all finite places $v$ of $L$ we have $\WD(R|_{\Gal(L_v^{ac}/L_v)}) = (\chi|_{L_v^\times} 
\circ \Art_{L_v}^{-1},0)$.
\end{itemize}
(See for instance \cite{alr}.) Then we see that $|\chi|^2=| \,\,\,|_L^{-k}$ (because both are characters
$\A_L^\times/(L_\infty^\times)^0 \ra \R^\times_{>0}$ which agree on $L^\times$) and so for every
finite place $v$ of $L$ the Weil-Deligne representation $\WD(R|_{\Gal(L_v^{ac}/L_v)})$ is
strictly pure of weight $k$. In particular $R$ is pure of weight $k$.

We have the following lemma.

\begin{lem}\label{gpure}
Suppose that $M/L$ is a finite extension of number fields. Suppose also that 
\[ R: \Gal(L^{ac}/L) \lra GL_n(\Q_l^{ac}) \]
is a continuous semisimple representation which is pure of weight $k$ and such that the restriction $R|_{\Gal
(L_x^{ac}/L_x)}$ is de Rham for all $x|l$. Suppose that 
\[ S: \Gal(M^{ac}/M) \lra GL_{an}(\Q_l^{ac}) \]
is another continuous representation with $S^\semis \cong
R|_{\Gal(M^{ac}/M)}^a$ for some $a \in \Z_{>0}$. Suppose finally that $w$ is a
place of $M$ above a finite place $v$ of $L$. If $\WD(S|_{\Gal(M_w^{ac}/M_w)})$
is pure of weight $k$, then $\WD(R|_{\Gal(L_v^{ac}/L_v)})$ is also pure of
weight $k$.
\end{lem}

\pfbegin
Write
\[ R|_{\Gal(M^{ac}/M)} = \bigoplus_i R_i \]
where each $R_i$ is irreducible. Then $\det R_i$ is de Rham at all places $x|l$ of $M$ and is
pure of weight $k \dim R_i$. Thus the top exterior power $\bigwedge^{\dim R_i}
\WD(R_i|_{\Gal(M_w^{ac}/M_w)})$ is
also pure of weight $k \dim R_i$. Lemma \ref{pure}(6) tells us that
\[ \WD(S|_{\Gal(M_w^{ac}/M_w)})^\Fsemis \cong \Big( \bigoplus_i
\WD(R_i|_{\Gal(M_w^{ac}/M_w)})^\Fsemis \Big)^a \cong 
\left( \WD(R|_{\Gal(M_w^{ac}/M_w)})^\Fsemis \right)^a, \]
and that $\WD(R|_{\Gal(M_w^{ac}/M_w)})^\Fsemis$ is pure of weight $k$.
Applying lemma \ref{pure}(1) and (2), we see that $\WD(R|_{\Gal(L_v^{ac}/L_v)})$
is also pure of weight $k$.
\pfend

\section{Shimura varieties}\label{s2}

In this section we recall some facts about the Shimura varieties
considered in \cite{ht} and prove proposition \ref{il2}.

In this section, 
\begin{itemize}
\item  let $E$ be an imaginary quadratic field, $F^+$ a totally real field and set $F=EF^+$;
\item let $p$ be a rational prime which splits as $p=uu^c$ in $E$;
\item and let $w=w_1,w_2,...,w_r$ be the places of $F$ above $u$;
\item and let $B$ be a division algebra with centre $F$ such that
\begin{itemize}
\item $\dim_F B=n^2$,
\item $B^\op \cong B \otimes_{F,c} F$,
\item at every place $x$ of $F$ either $B_x$ is split or a division algebra,
\item if $n$ is even then the number of finite places of $F^+$ above which
$B$ is ramified is congruent to $1+\frac{n}{2}[F^+:\Q]$ modulo $2$.
\end{itemize} \end{itemize}
Pick a positive involution $*$ on $B$ with $*|_F=c$. Let $V=B$ as a 
$B \otimes_F B^\op$-module. For $\beta \in B^{* = -1}$ define a pairing
\[ \begin{array}{rcl} ( \,\,\, ,\,\,\,): V \times V & \lra & \Q \\
(x_1,x_2) & \longmapsto & \tr_{F/\Q} \tr_{B/F} (x_1 \beta x_2^*). \end{array}\]
Also define an involution $\#$ on $B$ by $x^\#=\beta x^* \beta^{-1}$ and
a reductive group $G/\Q$ by setting, for any $\Q$-algebra $R$, the group
$G(R)$ equal to the set of 
\[ (\lambda , g) \in R^\times \times (B^\op \otimes_\Q R)^\times \]
such that
\[ gg^\# = \lambda. \]
Let $\nu: G \ra \G_m$ denote the multiplier character sending $(\lambda,g)$ to
$\lambda$. Note that if $x$ is a rational
prime which splits $x=yy^c$ in $E$ then
\[ \begin{array}{rcl} G(\Q_x) &\liso & (B_y^\op)^\times \times \Q_x^\times \\
(\lambda,g) &\longmapsto & (g_y,\lambda). \end{array} \]
We can and will assume that
\begin{itemize}
\item if $x$ is a rational prime which does not split in $E$ the $G \times
\Q_x$ is quasi-split;
\item the pairing $(\,\,\, ,\,\,\,)$ on $V \otimes_\Q \R$ has invariants
$(1,n-1)$ at one embedding $\tau:F^+ \into \R$ and invariants $(0,n)$ at
all other embeddings $F^+ \into \R$.
\end{itemize}
(See section I.7 of \cite{ht} for details.) 

Let $U$ be an open compact subgroup of $G(\A^\infty)$. Define a functor
$\gX_U$ from the category of pairs $(S,s)$, where $S$ is a connected locally
noetherian $F$-scheme and $s$ is a geometric point of $S$, to the category of sets, sending 
$(S,s)$ to the set of isogeny classes of four-tuples $(A,\lambda,i,\bareta)$
where
\begin{itemize}
\item $A/S$ is an abelian scheme of dimension $[F^+:\Q]n^2$;
\item $i:B \into \End(A) \otimes_\Z \Q$ such that $\Lie A \otimes_{(E \otimes_\Q
\CO_S), 1 \otimes 1} \CO_S$ is locally free over $\CO_S$ of rank $n$ and the two
actions of $F^+$ coincide;
\item $\lambda:A \ra A^\vee$ is a polarisation such that for all $b \in
B$ we have $\lambda \circ i(b) = i(b^*)^\vee \circ \lambda$;
\item $\bareta$ is a $\pi_1(S,s)$-invariant $U$-orbit of isomorphisms of
$B \otimes_\Q \A^\infty$-modules $\eta:V \otimes_\Q \A^\infty \ra VA_s$ which
take the standard pairing $(\,\,\, ,\,\,\,)$ on $V$ to a
$(\A^\infty)^\times$-multiple of the $\lambda$-Weil pairing on $VA_s$.
\end{itemize}
Here $VA_s = \Bigl(\invlim{N} A[N](k(s))\Bigr) \otimes_\Z \Q$ is the adelic Tate
module.
For the precise notion of isogeny class see section III.1 of \cite{ht}. If
$s$ and $s'$ are both geometric points of a connected locally noetherian
$F$-scheme $S$ then $\gX_U(S,s)$ and $\gX_U(S,s')$ are in canonical bijection.
Thus we may think of $\gX_U$ as a functor from connected locally noetherian
$F$-schemes to sets. We may further extend it to a functor from all
locally noetherian $F$-schemes to sets by setting
\[ \gX_U\biggl(\coprod_i S_i\biggr) = \prod_i \gX_U(S_i). \]
If $U$ is sufficiently small (i.e.\ for some finite place $x$ of $\Q$ the
projection of $U$ to $G(\Q_x)$ contains no element of finite order except
$1$) then $\gX_U$ is represented by a smooth projective variety $X_U/F$ of
dimension $n-1$. The inverse system of the $X_U$ for varying $U$ has a natural
right action of $G(\A^\infty)$.

We will write $\CA_{U}$ for the universal abelian variety over $X_{U}$. The action of
$G(\A^{\infty})$ on the inverse system of the $X_{U}$ extends to an action by
quasi-isogenies on the inverse system of the $\CA_{U}$. More precisely if $g^{-1}Vg
\subset U$ then we have a map $g:X_{V} \ra X_{U}$ and there is also a quasi-isogeny
of abelian varieties over $X_{V}$
\[ g:\CA_{V} \lra g^{*} \CA_{U}. \]
(By a quasi-isogeny from $A$ to $B$ we mean an element of $\Hom(A,B) \otimes_\Z
\Q$, some integer multiple of which is an isogeny.) 

Now suppose $l$ is a rational prime and $\xi$ is an irreducible algebraic representation of $G$ over $\Q_{l}^{ac}$. This defines a lisse $\Q_{l}^{ac}$-sheaf $\CL_{\xi} = \CL_{\xi,l}$ over
each $X_{U}$ and the action of $G(\A^{\infty})$ extends to these sheaves. (See section
III.2 of \cite{ht} for details.) We write
\[ H^{i}(X,\CL_{\xi,l}) = \dirlim{U} H^{i}(X_{U} \times_F F^{ac}, \CL_{{\xi,l}}). \]
It is a semi-simple admissible representation of $G(\A^{\infty})$ with a commuting continuous
action of $\Gal(F^{ac}/F)$. Thus we can write
\[ H^i(X,\CL_{\xi,l}) = \bigoplus_\pi \pi \otimes R^i_{\xi,l}(\pi) \]
where $\pi$ runs over irreducible admissible representations of $G(\A^\infty)$ over
$\Q_l^{ac}$, and $R^i_{\xi,l}(\pi)$ is a finite dimensional continuous representation of $\Gal(F^{ac}/F)$ over $\Q_l^{ac}$. We recall the following results from \cite{ht} (section VI.2 and corollary
V.6.2).

\begin{lem}\label{mo} Keep the above notation and assumptions. Suppose that $R_{\xi,l}^i(\pi)
\neq (0)$ and $R_{\xi,l}^j(\pi') \neq (0)$. If $\pi_x \cong \pi_x'$ for all but finitely
many finite places $x$ of $\Q$ then $\pi_x \cong \pi_x'$ for all finite places $x$ of 
$\Q$ which split in $E$. \end{lem}

\begin{lem}\label{coh} Keep the above notation and assumptions and choose $\imath:\Q_l^{ac}
\iso \C$. Suppose that $R^i_{\xi,l}(\pi)
\neq (0)$. Then there is an automorphic representation $\Pi$ of $GL_n(\A_F)$ which
occurs in the discrete spectrum and a character $\psi$ of $\A_E^\times/E^\times$
such that for all but finitely many finite places $x$ of $\Q$ we have $(\psi_x,\Pi_x)=
\BC(\imath \pi_x)$ in the notation of section VI.2 of \cite{ht}. (This characterises $(
\psi,\Pi)$ completely.)

If moreover $\Pi$ is cuspidal then the following hold.
\begin{enumerate}
\item $i=n-1$.
\item $\Pi^\vee = \Pi^c$.
\item $\Pi_\infty$ has the same infinitesimal character as an algebraic representation of 
the restriction of scalars from $F$ to $\Q$ of $GL_n$.
\item $\psi_\infty$ has the same infinitesimal character as an algebraic representation of 
the restriction of scalars from $E$ to $\Q$ of $GL_1$.
\item $\Pi_x$ is square integrable at some finite place of $F$. (In fact at any
place for which $B_x$ is not split.)
\item $R_{\xi,l}^{n-1}(\pi)^\semis \cong R_{l,\imath}(\Pi)^a \otimes R_{l,\imath} (\psi)$ for some $a \in \Z_{>0}$.
\end{enumerate} \end{lem}

We will write $(\psi,\Pi)=\BC_\imath(\pi)$. We will call $\BC_\imath(\pi)$ {\em cuspidal} 
if $\Pi$ is cuspidal.

For an irreducible algebraic representation $\xi$, in \cite{ht},
section III.2, integers $t_\xi, m_\xi \geq 0$ and an element
$\varepsilon_\xi\in \Q[S_{m_\xi}]$ (where $S_{m_\xi}$ is the symmetric
group of $m_\xi$ letters) are defined. We also set for each 
integer $N\geq 2$ set
\[ \varepsilon(m_\xi,N) = \prod_{x=1}^{m_\xi}\prod_{y\neq
  1}\frac{[N]_x-N^y}{N-N^y} \in \Q[(N^{\Z_{\geq 0}})^{m_\xi} ] , \]
 where $[N]_x$ denotes the endomorphism generated by multiplication by
$N$ on the $x$-th factor, and $y$ runs from $0$ to $2[F^+:\Q]n^2$ but
excluding 1. Also set 
\[ a_\xi = a_{\xi, N}= \varepsilon_\xi \varepsilon(m_\xi,N)^{2n-1} \in
\Q[(N^{\Z_{\geq 0}})^{m_\xi} \rtimes S_{m_\xi}] .\]
If we think of $(N^{\Z_{\geq 0}})^{m_\xi} \rtimes S_{m_\xi} \subset
\End(\CA_U^{m_\xi}/X_U)$ then $a_\xi \in \End(\CA_U^{m_\xi}/X_U)
\otimes_\Z \Q$. Write $\pi$ for the map $\CA_U^{m_\xi} \ra X_U$ and $\pi_i$ for the
composition of the $i^{th}$ inclusion $\CA_U \into \CA_U^{m_\xi}$ with $\pi$. Then
\begin{itemize}
\item $R^j\pi_* \Q_l^{ac} = \bigoplus_{i_1+...+i_{m_\xi}=j} \bigotimes_{a=1}^{m_\xi}
\wedge^{i_a} R^1\pi_{a \,\, *} \Q_l^{ac} $;
\item $\varepsilon(m_\xi,N)$ acts as an idempotent on $R^j\pi_* \Q_l^{ac}$;
\item $\varepsilon(m_\xi,N) R^j\pi_* \Q_l^{ac}= 
\begin{cases}
(0) & (j \neq m_\xi) \\ \bigotimes_{a=1}^{m_\xi} R^1\pi_{a \,\, *} \Q_l^{ac} & (j=m_\xi)
\end{cases}$;
\item $\varepsilon_\xi \varepsilon(m_\xi,N) R^{m_\xi}\pi_* \Q_l^{ac}(t_\xi) \cong \CL_\xi$;
\item $a_\xi$ acts as an idempotent on each $H^j(\CA_U^{m_\xi} \times_F F^{ac}, \Q_l^{ac}(t_\xi))$;
\item and $a_\xi H^j(\CA_U^{m_\xi} \times_F F^{ac},\Q_l^{ac}(t_\xi)) \cong 
\begin{cases}
(0) & (j< m_\xi)\\
H^{j-m_\xi}(X_U \times_F F^{ac}, \CL_\xi) & (j\geq m_\xi)
\end{cases}$.
\end{itemize}
Note that $a_\xi$ commutes with
the action of $G(\A^\infty)$. Moreover if $\imath:\Q_l^{ac} \cong \Q_{l'}^{ac}$ then
the same $a_\xi$ will work for $\xi$ and $\imath \xi$. 
(Note that in \cite{ht} the exponent $2n-1$ in the definition of
$a_\xi$ is erroneously missed. This gives
an idempotent on $R^j\pi_*\Q_l^{ac}(t_\xi)$ with the desired
properties, but we do not know that it is an idempotent on $H^j(\CA_U^{m_\xi} \times_F F^{ac}, \Q_l^{ac}(t_\xi))$. We only know that this group has a filtration of length
$2n-1$ (coming from the composite functor spectral sequence) on each of whose graded pieces $\varepsilon(m_\xi,N)$ is an idempotent.)

\begin{lem}\label{il3} Keep the above notation and assumptions. Let $\imath:
\Q_l^{ac} \iso \C$ and $\imath':\Q_{l'}^{ac}
\iso \C$. Suppose that $R^i_{\xi,l}(\pi) \neq (0)$ and that $\BC_\imath(\pi)$ is
cuspidal. Suppose also that $w$ is a place of $F$ above a split place of $E$ such that 
$B_w$ is split. If we assume conjecture \ref{il} then, writing $\xi'$ for $(\imath')^{-1} \imath \xi$ and $\pi'$ for $(\imath')^{-1}\imath \pi$,
\[ \imath \WD(R^{n-1}_{\xi,l}(\pi)|_{\Gal(F_w^{ac}/F_w)})^\semis \cong \imath' \WD(R^{n-1}_{
\xi',l'}(\pi')|_{\Gal(F_w^{ac}/F_w)})^\semis. \]
\end{lem}

\pfbegin Choose an open compact subgroup $U \subset G(\A^\infty)$ such that 
$\pi^U \neq (0)$. Also choose $e \in \Q_l^{ac}[U\backslash G(\A^\infty)/U]$ such that
$e$ is an idempotent on each $H^j(X_U \times_F F^{ac}, \CL_\xi)$ and
\[ e H^j(X_U \times_F F^{ac}, \CL_\xi) = \pi^U \otimes R^j_{\xi,l}(\pi), \]
and such that $e'=(\imath')^{-1} \imath e$ is an idempotent on each
$H^j(X_U \times_F F^{ac}, \CL_{\xi'})$ and
\[ e' H^j(X_U \times_F F^{ac},\CL_{\xi'}) = (\pi')^U \otimes R^j_ {\xi',l'}(\pi'). \]

Then, by conjecture \ref{il}, lemma \ref{coh} and the above discussion, for $\sigma \in W_{F_w}$ we have
\begin{align*}
& (\dim \pi^U) \tr\bigl(\sigma \big| \WD(R^{n-1}_{\xi,l}(\pi)|_{\Gal(F_w^{ac}/F_w)})\bigr) \\ 
=& \sum_j (-1)^{n-1-j} \tr \bigl(\sigma e \big| \WD\bigl(H^j(X_U \times_{F_w} F_w^{ac}, \CL_\xi)\bigr)\bigr) \\
=& \sum_j (-1)^{n-1+m_\xi -j} \tr \bigl(\sigma e a_\xi \big| \WD\bigl(H^j(\CA_U^{m_\xi} 
\times_{F_w} F_w^{ac}, \Q_l^{ac}(t_\xi))\bigr)\bigr) \\ 
=& \sum_j (-1)^{n-1+m_{\xi'} -j}
\tr \bigl(\sigma e' a_{\xi'} \big| \WD\bigl(H^j(\CA_U^{m_{\xi'}} \times_{F_w} F_w^{ac}, \Q_{l'}^{ac}(t_{\xi'}))\bigr)\bigr) \\
=& (\dim (\pi')^U) \tr\bigl(\sigma \big| \WD(R^{n-1}_{\xi',l'}(\pi')|_{\Gal(F_w^{ac}/F_w)})\bigr).
\end{align*}
The lemma follows.
\pfend

\begin{cor} Keep the above notation. Let $\imath: \Q_l^{ac} \iso \C$.
Suppose that $R^i_{\xi,l}(\pi) \neq (0)$ and that $\BC_\imath(\pi)=(\psi,\Pi)$ is
cuspidal. Suppose also that $w$ is a place of $F$ above a split place of $E$ such that 
$B_w$ is split. If we assume conjecture \ref{il} then 
\[ \imath \WD(R_{l,\imath}(\Pi)|_{\Gal(F_w^{ac}/F_w)})^\semis = \rec \bigl(\Pi_w^\vee |\det|^{\frac{1-n}{2}}\bigr)^\semis. \] \end{cor}

\pfbegin If $v \ndiv l$ this is part of theorem VII.1.9 of \cite{ht}. Thus suppose $v|l$. 
Choose $l' \neq l$ and $\imath':\Q_{l'}^{ac} \iso \C$. It suffices to show that
\[ \imath \WD(R_{l,\imath}(\Pi)|_{\Gal(F_w^{ac}/F_w)})^\semis \cong \imath' \WD(R_{l',\imath'}(\Pi)|_{\Gal(F_w^{ac}/F_w)})^\semis. \]
But this follows from lemmas \ref{coh} and \ref{il3}.
\pfend

We can now complete the proof of proposition \ref{il2}. Recall that $L$ is an
imaginary CM field and that $\Pi$ is a cuspidal automorphic representation
of $GL_n(\A_L)$ such that
\begin{itemize}
\item $\Pi \circ c \cong \Pi^\vee$;
\item $\Pi_\infty$ has the same infinitesimal character as some algebraic
representation over $\C$ of the restriction of scalars from $L$ to $\Q$ of
$GL_n$;
\item and for some finite place $x$ of $L$ the representation $\Pi_x$ is
square integrable.
\end{itemize}
Recall also that $v$ is a place of $L$ above a rational prime $p$. Choose a CM
field $L'$ which is a quadratic extension of $L$ in which $x$ and $v$ split. Also
choose primes $x'$ above $x$ and $v'$ above $v$ of $L'$, with $x' \ndiv v'(v')^c$. (This
is only important in the case that $x$ and $v$ lie above the same place of the maximal totally real subfield of $L$.) 

Choose an imaginary quadratic field $E$ not contained in $L'$, in which $p$ and the rational prime below $x$ split. Let $p=uu^c$ in $E$, let $F=EL'$ and let $F^+$ denote the maximal totally real subfield of $F$. Also choose places $z$ (resp.\ $w$) of $F$
above $x'$ (resp.\ $v'$) such that $w|_E=u$. 
Denote by $\Pi_F$ the base change of $\Pi$ to $GL_n(\A_F)$. Note that 
$\Pi_{F,w}$ is square integrable and hence $\Pi_F$ is cuspidal.

Choose a division algebra $B$ with centre $F$ as at the start of this section and
satisfying 
\begin{itemize}
\item $B_y$ is split for all places $y \ndiv  zz^c$ of $F$.
\end{itemize}
Also choose $\ast$, $\beta$ and $G$ as at the start of this section. Then it follows
from theorem VI.2.9 and lemma VI.2.10 of \cite{ht} that we can find
\begin{itemize}
\item a character $\psi: \A_E^\times/E^\times \ra \C^\times$,
\item an irreducible algebraic representation $\xi$ of $G$ over $\C$,
\item an irreducible representation $\pi$ of $G(\A^\infty)$ over $\C$,
\end{itemize}
such that
\begin{itemize}
\item $R^{n-1}_{\imath \xi}(\imath \pi) \neq (0)$,
\item $R^{n-1}_{\imath' \xi}(\imath' \pi) \neq (0)$,
\item and $\BC_\imath(\imath \pi)=\BC_{\imath'}(\imath' \pi) = (\psi,\Pi)$.
\end{itemize}
Proposition \ref{il2} now follows from the previous corollary.

\section{Integral models for Iwahori level structure}

We keep the notation of the last section. By a {\em geometric point} of a scheme $X$, we mean a map $s\colon \Spec k\ra X$, where $k$ is a field, such that the algebraic closure of the residue field of $s(\Spec k)$ is isomorphic to $k$. We call it a {\em closed geometric point} if its image $s(\Spec k)$ is a closed point of $X$.

Choose a maximal $\Z_{(p)}$-order $\CO_B$ of $B$
with $\CO_B^*=\CO_B$. Let $\CO_{F,w}$ be the integer ring of $F_w$, $\varpi_w$ a uniformiser for $\CO_{F,w}$, and $\CO_{B,w}=\CO_B\otimes_{\Z_{(p)}}\CO_{F,w}$. Also fix an isomorphism $\CO_{B,w}^\op \cong
M_n(\CO_{F,w})$, and let $\varepsilon \in B_w$ denote the element corresponding
to the diagonal matrix $(1,0,0,...,0) \in M_n(\CO_{F,w})$. We decompose
$G(\A^\infty)$ as
\[ G(\A^\infty) = G(\A^{\infty,p}) \times \biggl(\prod_{i=2}^r (B_{w_i}^\op)^\times\biggr) \times GL_n(F_w) \times \Q_p^\times. \]
For $m =(m_2,...,m_r) \in
\Z_{\geq 0}^{r-1}$, set
\[U^w_p(m)=\prod_{i=2}^r \ker \bigl((\CO_{B,w_i}^\op)^\times \ra (\CO_{B,w_i}^\op/w_i^{m_i})^\times\bigr) \subset \prod_{i=2}^r (B_{w_i}^\op)^\times.\]

Let $B_n$ denote the Borel subgroup of $GL_n$ consisting of upper triangular
matrices, and let $\Iw_{n,w}$ denote the subgroup
of $GL_n(\CO_{F,w})$ consisting of matrices which reduce modulo $w$ to
$B_n(k(w))$.
We will consider the following open subgroups of $G(\Q_p)$:
\begin{align*}
\Ma(m) &= U_p^w(m) \times GL_n(\CO_{F,w}) \times \Z_p^\times, \\
\Iw(m) &= U_p^w(m) \times \Iw_{n,w} \times \Z_p^\times.
\end{align*}
If $U^p$ is an open compact subgroup of $G(\A^{\infty,p})$, we will write $U_0$ (resp.\
$U$) for $U^p \times \Ma(m)$ (resp.\  $U^p \times \Iw(m)$).

We recall that in section III.4 \cite{ht} an integral model of $X_{U_0}$
over $\CO_{F,w}$ is defined (for $U^p$ sufficiently small). We denote the integral model also by $X_{U_0}$, by a slight abuse of notation (similarly for the integral model of $X_U$ defined later).
It represents a functor $\gX_{U_0}$
from locally noetherian $\CO_{F,w}$-schemes to sets. As above, $\gX_{U_0}$ is initially defined
on the category of connected locally noetherian $\CO_{F,w}$-schemes with a geometric point
to sets. It sends $(S,s)$ to the set of prime-to-$p$ isogeny classes of $(r+3)$-tuples 
$(A,\lambda,i,\bareta^p,\alpha_i)$, where
\begin{itemize}
\item $A/S$ is an abelian scheme of dimension $[F^+:\Q]n^2$;
\item $i:\CO_B \into \End(A) \otimes_\Z \Z_{(p)}$ such that $\Lie A
\otimes_{(\CO_{E,u} \otimes_{\Z_p} \CO_S), 1 \otimes 1} \CO_S$ is locally free
of rank $n$ and the two actions of $\CO_{F^+}$ coincide;
\item $\lambda: A \ra A^\vee$ is a prime-to-$p$ polarisation such that for all
$b \in \CO_B$ we have $\lambda \circ i(b)=i(b^*)^\vee \circ \lambda$;
\item $\bareta^p$ is a $\pi_1(S,s)$-invariant $U^p$-orbit of isomorphisms of
$B \otimes_\Q \A^{\infty,p}$-modules $\eta:V \otimes_\Q \A^{\infty,p} \ra V^pA_s$
which take the standard pairing $(\,\,\, ,\,\,\,)$ on $V$ to a
$(\A^{\infty,p})^\times$-multiple of the $\lambda$-Weil pairing on $V^pA_s$;
\item for $2 \leq i \leq r$, $\alpha_i:(w_i^{-m_i} \CO_{B,w_i}/
\CO_{B,w_i})_S \iso A[w_i^{m_i}]$ is an isomorphism of $S$-schemes with
$\CO_B$-actions;
\end{itemize}
Then $X_{U_0}$ is smooth and projective over $\CO_{F,w}$ (\cite{ht}, page 109).
As $U^p$ varies, the inverse system of the $X_{U_0}$'s
has an action of $G(\A^{\infty,p})$.

Given an $(r+3)$-tuple as above we will write $\CG_A$ for $\varepsilon A[w^\infty]$, a Barsotti-Tate
$\CO_{F,w}$-module. Over a base in which $p$ is nilpotent it is one dimensional and
{\em compatible}, i.e.\ the two actions of $\CO_{F,w}$ on $\Lie \CG_A$ coincide (see
\cite{ht}).  If $\CA_{U_0}$ denotes the universal abelian scheme over $X_{U_0}$,
we will write $\CG$ for $\CG_{\CA_{U_0}}$.

Write $\ol{X}_{U_0}$ for the special fibre $X_{U_0}\times_{\Spec \CO_{F,w}}\Spec k(w)$.
For $0\leq h\leq n-1$, we let $\ol{X}_{U_0}^{[h]}$ denote the reduced closed subscheme of $\ol{X}_{U_0}$
whose closed geometric points $s$ are those for which the maximal etale quotient of $\CG_s$ has
$\CO_{F,w}$-height at most $h$, and let 
\[\ol{X}_{U_0}^{(h)}=\ol{X}_{U_0}^{[h]}-\ol{X}_{U_0}^{[h-1]}\]
(where we set $\ol{X}_{U_0}^{[-1]}=\emptyset$). Then $\ol{X}_{U_0}^{(h)}$ is non-empty and smooth of
pure dimension $h$ (lemma III.4.3, corollary III.4.4 of \cite{ht}), and on it there is a short exact sequence
\[ (0) \lra \CG^0 \lra \CG \lra \CG^\et \lra (0) \]
where $\CG^0$ is a formal Barsotti-Tate $\CO_{F,w}$-module and $\CG^\et$ is an etale Barsotti-Tate
$\CO_{F,w}$-module with $\CO_{F,w}$-height $h$.

\begin{lem}\label{clos} If $0 \leq h \leq n-1$ then the Zariski closure of $\overline{X}^{(h)}_{U_0}$ contains $\overline{X}_{U_0}^{(0)}$. \end{lem}

\pfbegin This is `well known', but for lack of a reference we give a proof. Let $x$ be a closed geometric
point of $\overline{X}_{U_0}^{(0)}$. By lemma III.4.1 of \cite{ht} the formal completion of $\overline{X}_{U_0} \times \Spec k(w)^{ac}$ at $x$ is isomorphic to the equicharacteristic universal
deformation ring of $\CG_x$. According to the proof of proposition 4.2 of \cite{drin} this is isomorphic to 
$\Spf k(w)^{ac}[[T_1,...,T_{n-1}]]$ and we can choose the $T_i$ and a formal parameter $S$ on
the universal deformation of $\CG_x$ such that
\[ [\varpi_w](S) \equiv \varpi_wS+ \sum_{i=1}^{n-1} T_i S^{{\#k(w)}^i} + S^{{\#k(w)}^n}\ \ \bigl(\bmod\ S^{{\#k(w)}^n+1}\bigr).\]
Thus we get a morphism
\[ \Spec k(w)^{ac}[[T_1,...,T_{n-1}]] \lra \overline{X}_{U_0} \]
lying over $x\colon \Spec k(w)^{ac} \ra \overline{X}_{U_0}$, such that, if $k$ denotes the algebraic closure
of the field of fractions of $k(w)^{ac}[[T_1,...,T_{n-1}]]/(T_1,...,T_{n-h-1})$, then the induced map
\[ \Spec k \lra \overline{X}_{U_0} \]
factors through $\overline{X}_{U_0}^{(h)}$. Thus $x$ is in the closure of $\overline{X}_{U_0}^{(h)}$,
and the lemma follows. \pfend

Now we will define an integral model for $X_U$. It will represent a functor $\gX_U$
defined as follows.  Again we initially define it as a functor from the category of
connected locally noetherian schemes with a geometric point to sets, but then (as above) we extend
it to a functor from locally noetherian schemes to sets. The functor $\gX_U$ will send $(S,s)$ to the set of prime-to-$p$ isogeny
classes of $(r+4)$-tuples $(A,\lambda,i,\bareta^p, \CC, \alpha_i)$, where 
$(A,\lambda,i,\bareta^p, \alpha_i)$ is as in the definition of $\gX_{U_0}$ and
$\CC$ is a chain of isogenies 
\[\CC: \CG_A=\CG_0 \ra \CG_1 \ra \cdots \ra \CG_n = \CG_A/\CG_A[w]\]
of compatible Barsotti-Tate $\CO_{F,w}$-modules, each of degree $\#k(w)$ and with composite equal
to the canonical map $\CG_A \ra \CG_A/\CG_A[w]$. There is a natural transformation of functors $\gX_U \ra \gX_{U_0}$.

\begin{lem}\label{exist} If $U^p$ is sufficiently small, the functor $\gX_U$ is represented by
a scheme $X_U$ which is finite over $X_{U_0}$. The scheme
$X_U$ has some irreducible components of dimension $n$.
\end{lem}

\pfbegin 
By denoting the kernel of $\CG_0\ra \CG_j$ by $\CK_j\subset \CG[w]$, we can view the above chain as a flag
\[0=\CK_0\subset \CK_1\subset \CK_2\subset \cdots \subset \CK_{n-1}\subset \CK_n=\CG[w]\]
of closed finite flat subgroup schemes with $\CO_{F,w}$-action, with each $\CK_j/\CK_{j-1}$
having order $\#k(w)$.
Let $\CH$ denote the sheaf of Hopf algebras over $X_{U_0}$ defining $\CG[w]$. Then $\gX_U$ is
represented by a closed subscheme $X_U$ of the Grassmanian of chains of locally free direct
summands of $\CH$. (The closed conditions require that the subsheaves are sheaves of ideals defining a
flag of closed subgroup schemes with the desired properties.)  Thus $X_U$ is projective over
$\CO_{F,w}$. At each closed geometric point $s$ of $X_{U_0}$ the number of possible
$\CO_{F,w}$-submodules of $\CG[w]_s\cong \CG[w]^0_s\times \CG[w]^{\et}_s$ is finite, so
$X_U$ is finite over $X_{U_0}$. To see that $X_U$ has some components of dimension $n$ it suffices to note that on the generic fibre the map to $X_{U_0}$ is finite etale.
\pfend

We say an isogeny $\CG\ra \CG'$ of one-dimensional compatible Barsotti-Tate
$\CO_{F,w}$-modules of degree $\#k(w)$ over a scheme $S$ of characterstic $p$ has
{\em connected kernel} if it induces the zero map on $\Lie \CG$. We
will denote the Frobenius map by $F:\CG\ra \CG^{(p)}$ and let
$f=[k(w):\F_p]$, and then $F^f:\CG\ra \CG^{(\#k(w))}$ is an isogeny of
compatible Barsotti-Tate $\CO_{F,w}$-modules of degree $\#k(w)$ and has connected kernel. 

We have the following rigidity lemma.

\begin{lem}\label{rigid} Let $W$ denote the ring of integers of the completion of the maximal
unramified extension of $F_w$. Suppose that $R$ is an Artinian local $W$-algebra with residue field $k(w)^{ac}$. Suppose also that
\[ \CC: \CG_0 \ra \CG_1 \ra \cdots \ra \CG_n = \CG_0/\CG_0[w] \]
is a chain of isogenies of degree $\#k(w)$ of one-dimensional compatible formal Barsotti-Tate $\CO_{F,w}$-modules over $R$ of $\CO_{F,w}$-height $n$ with composite equal to multiplication by $\varpi_w$. If every
isogeny $\CG_{i-1} \ra \CG_i$ has connected kernel (for $i=1,...,n$) then $R$ is a $k(w)^{ac}$-algebra
and $\CC$ is the pull-back of a chain of Barsotti-Tate
$\CO_{F,w}$-modules over $k(w)^{ac}$, with all the isogenies isomorphic to $F^f$.
\end{lem}

\pfbegin
As the composite of the $n$ isogenies induces multiplication by $\varpi_w$ on the tangent space,
$\varpi_w=0$ in $R$, i.e.\ $R$ is a $k(w)^{ac}$-algebra. Choose a parameter $T_i$ for
$\CG_i$ over $R$. With respect to these choices, let $f_i(T_i)\in R[[T_i]]$ represent $\CG_{i-1} \ra
\CG_i$. We can write $f_i(T_i)=g_i(T_i^{p^{h_i}})$ with $h_i \in \Z_{\geq 0}$ and $g_i'(0) \neq 0$.
(See \cite{f}, chapter I, \S 3, Theorem 2.) As $\CG_{i-1} \ra \CG_i$ has connected kernel, $f_i'(0)=0$ and
$h_i>0$. As $f_i$ commutes with the action $[r]$ for all $r\in
\CO_{F,w}$, we have $\ol{r}^{p^{h_i}}=\ol{r}$ for all $\ol{r}\in
k(w)$, hence $h_i$ is a multiple of $f=[k(w):\F_p]$.  
Reducing modulo the maximal ideal of $R$ we see that $h_i \leq
f$ and so in fact $h_i=f$ and $g_i'(0) \in R^\times$. Thus $\CG_i
\cong \CG_0^{(\#k(w)^i)}$ in such a way that the isogeny $\CG_0 \ra
\CG_i$ is identified with $F^{fi}$. In particular
$\CG_0 \cong \CG_0^{(\#k(w)^n)}$
and hence $\CG_0 \cong \CG_0^{(\#k(w)^{nm})}$ for any $m \in \Z_{\geq 0}$. As $R$ is Artinian some power of
the absolute Frobenius on $R$ factors through $k(w)^{ac}$. Thus $\CG_0$ is a pull-back
from $k(w)^{ac}$ and the lemma follows.
\pfend

Now let $\overline{X}_U= X_U \times_{\Spec \CO_{F,w}} \Spec k(w)$ denote the special fibre of $X_U$, and let $Y_{U,i}$ denote the closed subscheme of $\ol{X}_U$
 over which $\CG_{i-1}\ra \CG_i$ has connected kernel.

The following proposition may be known to
experts. The case that $F_w/\Q_p$ is unramified follows from results of \cite{go},
however it is essential for our purposes to include the (harder) case where
$F_w/\Q_p$ is ramified. The proposition could be proved via the usual ``reduction to formal models''
argument and some linear algebra. However the following argument seems to be
easier and identifies formal parameters for the completion
of the strict Henselisation of $X_U$ at a geometric point of $X_U^{(0)}$. More precisely, 
they can be taken to be the scalars giving the linear maps $\Lie \CG_{i-1} \ra \Lie \CG_i$
with respect to some bases. 

\begin{prop} \label{semistable}
\begin{enumerate}
\item $X_U$ has pure dimension $n$ and semistable reduction over $\CO_{F,w}$, that is, for all closed
points $x$ of the special fibre $\ol{X}_U$, there exists an etale
morphism $V\ra X_U$ with $x$ contained in the image of $V$ and an etale $\CO_{F,w}$-morphism:
\[ V\lra \Spec \CO_{F,w}[T_1,...,T_n]/(T_1\cdots T_m-\varpi_w) \]
for some $1\leq m\leq n$, where $\varpi_w$ is a uniformizer of $\CO_{F,w}$.
\item $X_U$ is regular and the natural map $X_U\ra X_{U_0}$ is finite and flat.
\item Each $Y_{U,i}$ is smooth over $\Spec k(w)$ of pure dimension $n-1$,
$\ol{X}_U=\bigcup_{i=1}^nY_{U,i}$ and, for $i \neq j$ the schemes $Y_{U,i}$
and $Y_{U,j}$ share no common connected component. In particular, $X_U$ has
strictly semistable reduction.
\end{enumerate}
\end{prop}

\pfbegin 
In this proof we will make repeated use of the following version of Deligne's
homogeneity principle (\cite{dr}). Write $W$ for the ring of integers of the
completion of the maximal unramified extension of $F_w$. In what follows, if
$s$ is a closed geometric point of an $\CO_{F,w}$-scheme $X$ locally of finite
type, then we write $\CO_{X,s}^\wedge$ for the completion of the strict
Henselisation of $X$ at $s$, i.e.\ $\CO_{X\times\Spec W,s}^\wedge$. Let
$\bf P$ be a
property of complete noetherian local $W$-algebras such that if $X$ is a
$\CO_{F,w}$-scheme locally of finite type then the set of closed geometric
points $s$ of $X$ for which $\CO_{X,s}^\wedge$ has property $\bf P$
is Zariski open. Also let $X\ra X_{U_0}$ be a finite morphism with the
following properties
\begin{enumerate}
\renewcommand{\labelenumi}{(\roman{enumi})}
\item If $s$ is a closed geometric point of $\ol{X}_{U_0}^{(h)}$ then, up to
isomorphism, $\CO_{X,s}^\wedge$ does not depend on $s$ (but only on $h$).
\item There is a unique geometric point of $X$ above any geometric point of 
$\ol{X}_{U_0}^{(0)}$.
\end{enumerate}
If $\CO_{X,s}^\wedge$ has property $\bf P$ for every geometric point
of $X$ over $\ol{X}_{U_0}^{(0)}$, then $\CO_{X,s}^\wedge$ has
property $\bf P$ for every closed geometric point of $X$. Indeed, if we let $Z$ denote
the closed subset of $X$ where {\bf P} does not hold, then its image in
$X_{U_0}$ is closed and either is empty or contains some $\ol{X}_{U_0}^{(h)}$.
In the latter case, by lemma \ref{clos}, it also contains $\ol{X}_{U_0}^{(0)}$, which is impossible.
Thus $Z$ must be empty.

Note that both $X=X_U$ and $X=Y_{U,i}$ satisfy the above condition (ii) for
the homogeneity principle, by letting $R=k(w)^{ac}$ in lemma \ref{rigid}.

(1): The dimension of $\CO_{X_U,s}^\wedge$ as $s$ runs over geometric points of $X_U$ above
$\ol{X}_{U_0}^{(0)}$ is constant, say $m$. Applying the homogeneity principle to $X=X_U$ with
$\bf P$ being `dimension $m$', we see that $X_U$ has pure dimension $m$. By lemma \ref{exist} we
must have $m=n$ and $X_U$ has pure dimension $n$.

Now we will apply the above homogeneity principle to $X=X_U$ taking $\bf P$
to be `isomorphic to $W[[T_1,...,T_n]]/(T_1\cdots T_m-\varpi_w)$ for some
$m \leq n$'. By a standard argument (see e.g.\ the proof of proposition 4.10
of \cite{y}) the set of points with this property is open and if all closed
geometric points of $X_U$ have this property $\bf P$ then $X_U$ is semistable of pure dimension
$n$.

Let $s$ be a geometric point of $X_U$ over a point of $\ol{X}_{U_0}^{(0)}$.
Choose a basis $e_i$ of $\Lie \CG_i$ over $\CO_{X_U,s}^\wedge$ such
that $e_n$ maps to $e_0$ under the isomorphism $\CG_n =\CG_0/\CG_0[w]
\iso \CG_0$ induced by $\varpi_w$. With
respect to these bases let $X_i \in \CO_{X_U,s}^\wedge$ represent the linear map $\Lie \CG_{i-1} \ra \Lie \CG_i$. Then
\[ X_1\cdots X_n = \varpi_w. \]
Moreover it follows from lemma \ref{rigid} that $\CO_{X_U,s}^\wedge/(X_1,...,X_n) = k(w)^{ac}$.
(Because, by lemma III.4.1 of \cite{ht}, $\CO_{X_{U_0},s}^\wedge$ is the
universal deformation space of $\CG_s$. Hence by lemma \ref{rigid}, $\CO_{X_U,s}^\wedge$ is the universal deformation space of the chain
\[ \CG_s \stackrel{F^f}{\lra} \CG_s^{(\#k(w))} \stackrel{F^f}{\lra} \cdots 
\stackrel{F^f}{\lra} \CG_s^{(\#k(w)^n)} \cong \CG_s/\CG_s[\varpi_w].) \]
Thus we get a surjection
\[ W[[X_1,...,X_n]]/(X_1\cdots X_n-\varpi_w) \onto \CO_{X_U,s}^\wedge \]
and as $\CO_{X_U,s}^\wedge$ has dimension $n$ this map must be an isomorphism.

(2): We see at once that $X_U$ is regular. Then \cite{ak} V, 3.6 tells us that $X_U \ra X_{U_0}$ is
flat.

(3): We apply the homogeneity principle to $X=Y_{U,i}$ taking $\bf P$ to be
`formally smooth of dimension $n-1$'. If $s$ is a geometric point of $Y_{U,i}$
above $\ol{X}_{U_0}^{(0)}$ then we see that $\CO_{Y_{U,i},s}^\wedge$
is cut out in $\CO_{X_U,s}^\wedge \cong W[[X_1,...,X_n]]/(X_1\cdots X_n-\varpi_w)$ by the
single equation $X_i=0$. (We are using the parameters $X_i$ defined above.) Thus
\[ \CO_{Y_{U,i},s}^\wedge \cong k(w)^{ac}[[X_1,...,X_{i-1},X_{i+1},...,X_n]] \]
is formally smooth of dimension $n-1$. We deduce that $Y_{U,i}$ is smooth of
pure dimension $n-1$. 

As our $\CG/\ol{X}_U$ is one-dimensional, over a closed point, at least one of the isogenies
$\CG_{i-1}\ra \CG_i$ must have connected kernel, which shows that $\ol{X}_U=\bigcup_iY_{U,i}$.
Suppose $Y_{U,i}$ and $Y_{U,j}$ share a connected component $Y$ for some $i\neq j$. Then
$Y$ would be finite flat over $\ol{X}_{U_0}$ and so the image of $Y$ would meet 
$\ol{X}_{U_0}^{(n-1)}$. This is impossible, because above a closed point of $\ol{X}_{U_0}^{(n-1)}$ only one
isogeny among the chain can have connected kernel. Thus, for $i \neq j$ the closed subschemes $Y_{U,i}$ and $Y_{U,j}$ have no connected component in common.
\pfend

By the strict semistability, if we write, for $S \subset \{1,...,n\}$,
\[ Y_{U,S} = \bigcap_{i \in S} Y_{U,i},\ \ Y_{U,S}^0 =
Y_{U,S}-\bigcup_{T \supsetneq S} Y_{U,T}\]
then $Y_{U,S}$ is smooth over $\Spec k(w)$ of pure dimension $n-\# S$ and
$Y_{U,S}^0$ are disjoint for different $S$. With respect to the finite flat
map $\ol{X}_U\ra \ol{X}_{U_0}$, the inverse image of $\ol{X}_{U_0}^{[h]}$ is
exactly the locus where at least $n-h$ of the isogenies have connected kernel,
i.e.\ $\bigcup_{\# S\geq n-h}Y_{U,S}$. Hence the inverse image of
$\ol{X}_{U_0}^{(h)}$ is equal to $\bigcup_{\# S=n-h}Y_{U,S}^0$.

The inverse systems $X_U \ra X_{U_0}$, as $U^p$ varies, have compatible actions of 
$G(\A^{\infty,p})$. For any $S$, the systems of subvarieties $Y_{U,S}$ and $Y_{U,S}^0$
are stable under this action. As in characterstic zero, these actions extend to
actions on the universal abelain varieties $\CA_{U_0}$ and $\CA_U$ over these bases (which we again denote by the same symbols $\CA_{U_0}$ and $\CA_U$).
This action is by prime-to-$p$ quasi-isogenies. (By a prime-to-$p$ quasi-isogeny we mean a quasi-isogeny such that the multiple by some integer prime to $p$ is an isogeny of degree prime to $p$.) 

Let $l$ be a prime and $\xi$ be an irreducible representation of $G$ over $\Q_l^{ac}$.
If $l \neq p$ then the sheaf $\CL_\xi$ extends to a lisse sheaf on our integral models
of $X_{U_0}$ and $X_U$ (using exactly the same construction as in characteristic
zero), and $a_\xi=a_{\xi,N}\in \End(\CA_U^{m_\xi}/X_U)\otimes_\Z \Q$ extends to the $\CA_U^{m_\xi}$ over integral models. We take $N$ prime to $p$, so that $a_\xi$ extends as etale morphisms on $\CA_U^{m_\xi}$. Also we will denote $\CA^{m_\xi}_U\times_{X_U}Y_{U,S}$ by $\CA^{m_\xi}_{U,S}$ for simplicity. We make the following definitions. 
\begin{itemize}
\item Define the admissible $G(\A^{\infty,p})$-modules with a commuting
continuous action of $\Gal(F^{ac}/F)$:
\begin{align*}
H^j(X_{\Iw(m)}, \CL_\xi) &= \dirlim{U^p} H^j(X_U \times_{F} F^{ac}, \CL_\xi) =
H^j(X, \CL_\xi)^{\Iw(m)}\\
H^j(\CA^{m_\xi}_{\Iw(m)}, \Q_l^{ac}) &= \dirlim{U^p} H^j(\CA_U^{m_\xi} \times_{F} 
F^{ac}, \Q_l^{ac}).
\end{align*}
\item If $l \neq p$, define admissible
$G(\A^{\infty,p}) \times \Frob_w^\Z$-modules:
\begin{align*}
H^j(Y_{\Iw(m),S}, \CL_\xi) &= \dirlim{U^p} H^j(Y_{U,S} 
\times_{k(w)} k(w)^{ac}, \CL_\xi)\\
H^j_c(Y_{\Iw(m),S}^0, \CL_\xi) &= \dirlim{U^p} H^j_c(Y_{U,S}^0 \times_{k(w)} k(w)^{ac}, \CL_\xi)\\
H^j(\CA^{m_\xi}_{\Iw(m),S},\Q_l^{ac}) &= \dirlim{U^p} H^j(\CA^{m_\xi}_{U,S} \times_{k(w)} k(w)^{ac}, \Q_l^{ac}).
\end{align*}
\item If $l=p$ and $\tau:W_0 \into \Q_l^{ac}$ over $\Z_p=\Z_l$, set (let $W_0$ denote the Witt ring of $k(w)$):
\[ H^j(\CA^{m_\xi}_{\Iw(m),S}/W_0) \otimes_{W_0, \tau} \Q_l^{ac}
= \dirlim{U^p} H^j(\CA^{m_\xi}_{U,S}/W_0) \otimes_{W_0, \tau}
\Q_l^{ac},\]
an admissible $G(\A^{\infty,p}) \times \Frob_w^\Z$-module. (Here
$ H^j(\CA^{m_\xi}_{U,S}/W_0)$ denote crystalline cohomology and we
let $\Frob_w$ act by the $[k(w):\F_p]$-power of the crystalline Frobenius.)
\end{itemize}
Note that if $l \neq p$ then $a_\xi$ is an idempotent on
$H^j(\CA^{m_\xi}_{\Iw(m),S},\Q_l^{ac})$ and  $a_xi
H^j(\CA^{m_\xi}_{\Iw(m),S},\Q_l^{ac}) = H^j(Y_{\Iw(m),S}, \CL_\xi)$ (for the
same reason this is true in characteristic zero). Similarly 
$a_\xi$ defines an idempotent on each
$H^j(\CA^{m_\xi}_{U,S}/W_0)$ and hence on
$H^j(\CA^{m_\xi}_{\Iw(m),S}/W_0) \otimes_{W_0, \tau} \Q_l^{ac}$, by the crystalline analogue of the
same argument.

We will call two irreducible admissible representations $\pi$ and $\pi'$
of $G(\A^{\infty,p})$ {\em nearly equivalent} if $\pi_x \cong \pi_x'$ for
all but finitely many rational primes $x \ndiv p$. If $M$ is an admissible
$G(\A^{\infty, p})$-module and $\pi$ is an irreducible admissible representation
of $G(\A^{\infty, p})$ then we define the $\pi$-{\em near isotypic} component
$M[\pi]$ of $M$ to be the largest $G(\A^{\infty, p})$-submodule of $M$ all
whose irreducible subquotients are nearly equivalent to $\pi$. Then
\[ M = \bigoplus M[\pi] \]
as $\pi$ runs over representatives of near equivalence classes of irreducible admissible
$G(\A^\infty)$-modules. (This follows from the following fact. Suppose that
$A$ is a (commutative) polynomial algebra over $\C$ in countably many
variables, and that $M$ is an $A$-module which is finitely generated over
$\C$. Then we can write 
\[ M= \bigoplus_\gm M_\gm, \]
where $\gm$ runs over maximal ideals of $A$ with residue field $\C$.)

Note that as $\CA_U$ is smooth over $X_U$ the varieties $\CA_U^{m_\xi}$ are strictly semistable. The special fibre is the union of the smooth
subschemes $\CA_{U,i}^{m_\xi}=\CA_U^{m_\xi} \times_{X_U} Y_{U,i}$ for $1\leq i\leq n$, and for $i \neq j$ the subschemes
$\CA_{U,i}^{m_\xi}$ and $\CA_{U,j}^{m_\xi}$ have no component in
common. Now we have:

\begin{prop} \label{specseq}
Suppose that $\pi$ is an irrreducible admissible representation of $G(\A^{\infty,p})$. For each rational prime $l$, there is a spectral sequence
\[ E_1^{i,j}(\Iw(m),\xi)[\pi] \Rightarrow \WD\bigl(H^{i+j}(X_{\Iw(m)}, \CL_\xi)|_{\Gal(F_w^{ac}/F_w)}\bigr)[\pi] \]
where $\displaystyle E_1^{i,j}(\Iw(m),\xi)=\bigoplus_{s \geq \max (0,-i)} \bigoplus_{\# S = i+2s+1}H^j_S$ and
\[ H^j_S=
\begin{cases}
a_\xi H^{j+m_\xi-2s}(\CA^{m_\xi}_{\Iw(m),S}, \Q_l^{ac}(t_\xi-s))=H^{j-2s}(Y_{\Iw(m),S}, \CL_\xi(-s)) & (l\neq p)\\
a_\xi H^{j+m_\xi-2s}(\CA^{m_\xi}_{\Iw(m),S}/W_0) \otimes_{W_0,\tau} 
\Q_l^{ac}(t_\xi-s) & (l=p)
\end{cases}. \]
\end{prop}

\pfbegin
We will use the functoriality of weight spectral sequences with
respect to the pull back by etale morphisms. This follows immediately
from the etale local nature of the construction of weight spectral
sequences. For $\ell=p$, this property has been used in \cite{saitomf}
and \cite{och}. For $\ell\neq p$ this is proven for general morphisms
(not necessarily etale) in \cite{saito}.

If $l \neq p$ then the Rapoport-Zink weight spectral sequence (\cite{rz}, \cite{saito}) is a spectral sequence
\begin{align*}
E_1^{i,j}= & \bigoplus_{s \geq \max (0,-i)} \bigoplus_{\# S = i+2s+1}
H^{j-2s}(\CA^{m_\xi}_{U,S}\times_{k(w)} k(w)^{ac}, \Q_l^{ac}(t-s)) \\
&\Rightarrow H^{i+j}(\CA_U^{m_\xi} \times_F F_w^{ac}, \Q_l^{ac}(t)).
\end{align*}
Taking $t=t_\xi$, applying $a_\xi$, replacing $j$ by $j+m_\xi$,
and passing to the direct limit over $U^p$ we get a spectral sequence of 
$G(\A^{\infty,p})$-modules
\begin{align*}
E_1^{i,j}(\Iw(m),\xi) = &\bigoplus_{s \geq \max (0,-i)} \bigoplus_{\# S = i+2s+1}
a_\xi H^{j+m_\xi-2s}(\CA^{m_\xi}_{\Iw(m),S}, \Q_l^{ac}(t_\xi-s)) \\
& \Rightarrow a_\xi H^{i+j+m_\xi}(\CA^{m_\xi}_{\Iw(m)}, \Q_l^{ac}(t_\xi)) = H^{i+j}(X_{\Iw(m)}, \CL_\xi),
\end{align*}
or equivalently
\[ E_1^{i,j}(\Iw(m),\xi) = \bigoplus_{s \geq \max (0,-i)} \bigoplus_{\# S = i+2s+1}
H^{j-2s}(Y_{\Iw(m),S}, \CL_\xi(-s)) \Rightarrow H^{i+j}(X_{\Iw(m)}, \CL_\xi).  \]
Hence we get the desired spectral sequence after passing to $\pi$-near isotypic components and identifying 
$a_\xi H^{i+j+m_\xi}(\CA^{m_\xi}_{\Iw(m)}, \Q_l^{ac}(t_\xi))$
(resp.\  $H^{i+j}(X_{\Iw(m)}, \CL_\xi)$) with their associated Weil-Deligne representations. Note that $I_{F_w}$ acts trivially on these spaces, the spectral sequence is 
equivariant for the action of $\Frob_w^\Z$ and the endomorphism $N$ is induced by the identity map 
\begin{align*}
N \colon & \bigoplus_{\# S = i+2s+1}
a_\xi H^{j+m_\xi-2s}(\CA^{m_\xi}_{\Iw(m),S}, \Q_l^{ac}(t_\xi-s)) \liso  \\
 & \bigoplus_{\# S = (i+2)+2(s-1)+1} a_\xi H^{(j-2)+m_\xi-2(s-1)}(
\CA^{m_\xi}_{\Iw(m),S}, \Q_l^{ac}(t_\xi-(s-1)))
\end{align*}
(resp.\
\[ N \colon \bigoplus_{\# S = i+2s+1} H^{j-2s}(Y_{\Iw(m),S}, \CL_\xi(-s)) \liso \!\!\!\!\!\!\!\! \bigoplus_{\# S = (i+2)+2(s-1)+1} \!\!\!\!\!\!\!\! H^{(j-2)-2(s-1)}(Y_{\Iw(m),S}, \CL_\xi(1-s))). \] 

If $l=p$ then the Mokrane \cite{mok} weight spectral sequence is a spectral
sequence
\[ E_1^{i,j} = \bigoplus_{s \geq \max (0,-i)} \bigoplus_{\# S = i+2s+1}
H^{j-2s}(\CA^{m_\xi}_{U,S}/W_0)(-s) \Rightarrow H^{i+j}((\CA_U^{m_\xi})^\times/
W_0^\times),  \]
computing the log-crystalline cohomology of $\CA_U^{m_\xi}$ in terms of the crystalline
cohomology of the $\CA^{m_\xi}_{U,S}$. Combining this with Tsuji's
comparison theorem \cite{tsuji} we get, for any choice of an embedding $\tau:
W_0 \into \Q_l^{ac}$ over $\Z_p=\Z_l$, a spectral sequence
\begin{align*}
E_1^{i,j} = & \bigoplus_{s \geq \max (0,-i)} \bigoplus_{\# S = i+2s+1}
H^{j-2s}(\CA^{m_\xi}_{U,S}/W_0) \otimes_{W_0, \tau} \Q_l^{ac}(t-s) \\
& \Rightarrow  \WD\bigl(H^{i+j}(\CA_U^{m_\xi} \times_F F_w^{ac}, \Q_l^{ac}(t))\bigr). 
\end{align*}
Taking $t=t_\xi$, applying $a_\xi$ (which is a linear combination of
etale morphisms by our choice of $N$), replacing $j$ by $j+m_\xi$,
and passing to the direct limit over $U^p$ we get a spectral sequence of 
$G(\A^{\infty,p})$-modules
\begin{align*}  
E_1^{i,j}(\Iw(m),\xi) = & \bigoplus_{s \geq \max (0,-i)} \bigoplus_{\# S = i+2s+1}
a_\xi H^{j+m_\xi-2s}(\CA^{m_\xi}_{\Iw(m),S}/W_0) \otimes_{W_0,\tau} 
\Q_l^{ac}(t_\xi-s) \\ 
& \Rightarrow \WD(H^{i+j}(X_{\Iw(m)}, \CL_\xi)|_{\Gal(F_w^{ac}/F_w)}).
\end{align*}
Hence we get the desired spectral sequence after passing to $\pi$-near
isotypic components (because $G(\A^{\infty,p})$ acts by etale morphisms).
On $\WD(H^{i+j}(X_{\Iw(m)}, \CL_\xi)|_{\Gal(F_w^{ac}/F_w)})$, the inertia group $I_{F_w}$ acts trivially, the action of $\Frob_w$ is compatible with the action of the crystalline Frobenius on $a_\xi H^{j-2s}(\CA^{m_\xi}_{\Iw(m),S}/W_0)(t_\xi-s)$, and
the endomorphism $N$ of $\WD(H^{i+j}(X_{\Iw(m)}, \CL_\xi)|_{\Gal(F_w^{ac}/F_w)})$ is induced by the
identity maps
\begin{align*} N\colon &\bigoplus_{\# S = i+2s+1}
a_\xi H^{j+m_\xi-2s}(\CA^{m_\xi}_{\Iw(m),S}/W_0)(t_\xi-s) \otimes_{W_0, \tau} \Q_l^{ac} \liso \\
&\bigoplus_{\# S = (i+2)+2(s-1)+1} a_\xi H^{(j-2)+m_\xi-2(s-1)}
(\CA^{m_\xi}_{\Iw(m),S}/W_0)(t_\xi-(s-1)) \otimes_{W_0, \tau}{\Q_l^{ac}}.
\end{align*}
\pfend

\section{Computing the cohomology of $Y_{U,S}$.}

In  this section we use the results of \cite{ht} to compute $H^j(Y_{\Iw(m),S},\CL_\xi)$.
As a result we can show that a large part of the above spectral sequences 
degenerate at $E_1$ and from this we deduce our main theorems.

We will keep the notation of the last section.

First we will relate the open strata $Y_{U,S}^0$ to the Igusa varieties of the
first kind defined in \cite{ht}. For $0\leq h\leq n-1$, $m_1\in \Z_{\geq 0}$ and $m\in \Z^{r-1}_{\geq 0}$, we write $I^{(h)}_{U^p,(m_1,m)}$ for the Igusa varities of the
first kind defined on page 121 of \cite{ht}. We also define an {\it
Iwahori-Igusa variety of the first kind}
\[I^{(h)}_U/\ol{X}_{U_0}^{(h)}\]
as the moduli space of chains of isogenies 
\[ \CG^{\et} =\CG_0\ra \CG_1 \ra \cdots \ra \CG_h=\CG^{\et}/\CG^{\et}[w]\]
of etale Barsotti-Tate $\CO_{F,w}$-modules, each isogeny having degree
$\#k(w)$ and with composite equal to the natural map $\CG^\et \ra \CG^\et/
\CG^\et[w]$. Then $I^{(h)}_U$ is finite etale over $\ol{X}_{U_0}^{(h)}$, and
as the Igusa variety $I^{(h)}_{U^p,(1,m)}$ classifies the isomorphisms
\[\alpha_1^{\et}:(w^{-1}\CO_{F,w}/\CO_{F,w})^h_{\ol{X}_{U_0}^{(h)}}\lra\CG^{\et}[w],\]
the natural map
\[ I^{(h)}_{U^p,(1,m)} \lra I^{(h)}_U \]
is finite etale and Galois with Galois group $B_h(k(w))$. Hence we can
identify $I^{(h)}_U$ with $I^{(h)}_{U^p,(1,m)}/B_h(k(w))$. Note that the
system $I^{(h)}_U$ for varying $U^p$ naturally inherits the action of
$G(\A^{\infty,p})$.

\begin{lem} \label{finbij}
For $S\subset \{1,..., n\}$ with $\#S=n-h$, there exists a finite map of
$\ol{X}^{(h)}_{U_0}$-schemes
\[ \varphi:Y_{U,S}^0 \lra I^{(h)}_U\]
which is bijective on the geometric points.
\end{lem}

\pfbegin
The map is defined in a natural way from the chain of isogenies $\CC$ by
passing to the etale quotient $\CG^{\et}$, and it is finite as $Y_{U,S}^0$
(resp.\ $I^{(h)}_U$) is finite (resp.\ finite etale) over $\ol{X}^{(h)}_{U_0}$.
Let $s$ be a closed geometric point of $I^{(h)}_U$ with a chain  of isogenies
\[ \CG_s^\et =\CG^\et_0 \ra \cdots \ra \CG^\et_h=\CG_s^\et/\CG_s^\et[w]. \]
For $1 \leq i \leq n$ let $j(i)$ denote the number of elements of $S$ which
are less than or equal to $i$. Set $\CG_i=(\CG_s^0)^{(\#k(w)^{j(i)})} \times
\CG^\et_{i-j(i)}$. If $i \not\in S$, define an isogeny $\CG_{i-1} \ra
\CG_i$ to be the identity times the given isogeny $\CG_{i-1}^\et \ra
\CG^\et_i$. If $i \in S$, define an isogeny $\CG_{i-1} \ra \CG_i$ to be
$F^f$ times the identity. Then
\[ \CG_0 \ra \cdots \ra \CG_n \]
defines the unique geometric point of $Y_{U,S}^0$ above $s$.
\pfend

\begin{cor} \label{geom}
Suppose that $l \neq p$. 
For every $S\subset \{1,..., n\}$ with $\#S=n-h$ and every $i\in \Z_{\geq 0}$, we have isomorphisms
\begin{align*}
H^i_c(Y_{U,S}^0\times_{k(w)}k(w)^{ac},\CL_\xi) &\liso H^i_c(I^{(h)}_U\times_{k(w)}k(w)^{ac},\CL_\xi)\\
&\liso H^i_c(I^{(h)}_{U^p,(1,m)}\times_{k(w)}k(w)^{ac},\CL_\xi)^{B_h(k(w))}
\end{align*}
that are compatible with the actions of $G(\A^{\infty,p})$ when we vary $U^p$.
\end{cor}

\pfbegin By lemma \ref{finbij}, for any lisse $\Q_l^{ac}$-sheaf $\CF$ on
$I^{(h)}_U$, we have $\CF\cong \varphi_*\varphi^*\CF$ by looking at the stalks
at all geometric points. As $\varphi$ is finite the first isomorphism follows.
The second isomorphism follows easilly as $I^{(h)}_{U^p,(1,m)} \ra
I^{(h)}_U$ is finite etale and Galois with Galois group $B_h(k(w))$.
\pfend

If $l \neq p$, we define
\[ H_c^j(I_{\Iw(m)}^{(h)},\CL_\xi) = \dirlim{U^p} H^j_c(I^{(h)}_U \times_{k(w)} k(w)^{ac}, \CL_\xi) = H_c^j(I^{(h)},\CL_\xi)^{U^w_p \times \Iw_{h,w}} \]
in the notation of page 136 of \cite{ht}. It is an admissible $G(\A^{\infty,p}) \times
\Frob_w^\Z$-module. In the notation of page 136 of \cite{ht}, $\Frob_w$ acts as
\[(1,p^{-[k(w):\F_p]},-1,1,1)\in G(\A^{\infty,p})\times (\Q_p^\times/
\Z_p^\times) \times \Z\times GL_h(F_w)\times \biggl(\prod_{i=2}^r
(B_{w_i}^\op)^\times\biggr),\]
where we have identified $D_{F_w,n-h}^\times/\CO_{D_{F_w,n-h}}^\times$ with
$\Z$ via $w\circ \det$.
We also define elements of $\Groth(G(\A^{\infty,p})
\times \Frob_w^\Z)$ (we write $\Groth(G)$ for the Grothendieck group of
admissible $G$-modules) as follows:
\begin{align*}
\bigl[H(Y_{\Iw(m),S},\CL_\xi)\bigr] &= \sum_i (-1)^{n-\# S -i}  H^i(Y_{\Iw(m),S},\CL_\xi),\\
\bigl[H_c(Y_{\Iw(m),S}^0,\CL_\xi)\bigr] &= \sum_i (-1)^{n-\# S -i}  H^i_c(Y_{\Iw(m),S}^0,
\CL_\xi),\\
\bigl[H_c(I_{\Iw(m)}^{(h)},\CL_\xi)\bigr] &= \sum_i (-1)^{h -i}  H^i_c(I_{\Iw(m)}^{(h)}, \CL_\xi).
\end{align*}
Finally set
\[ [H(X,\CL_\xi)] = \sum_j (-1)^{n-1-j} H^j(X,\CL_\xi) \in \Groth(G(\A^\infty)). \]

Theorem V.5.4 of \cite{ht} tells us that (for $l \neq p$)
\begin{align*}
n\bigl[H_c(I_{\Iw(m)}^{(h)},\CL_\xi)\bigr] &= 
n\bigl[H_c(I^{(h)},\CL_\xi)\bigr]^{U^w_p(m)\times \Iw_{h,w}}\\
&= \sum_i (-1)^{n-1-i} \Red^{(h)}\bigl[H^i(X,\CL_\xi)^{U^w_p(m)}\bigr]
\end{align*}
in $\Groth(G(\A^{\infty,p}) \times \Frob_w^\Z)$, where
\[ \Red^{(h)}: \Groth(GL_n(F_w) \times \Q_p^\times) \lra \Groth(\Frob_w^\Z) \]
is the composite of the normalised Jacquet functor
\[ J_{N_h^\op}: \Groth(GL_n(F_w) \times \Q_p^\times) \lra \Groth(GL_{n-h}(F_w)
\times GL_h(F_w) \times \Q_p^\times) \]
with the functor
\[ \Groth(GL_{n-h}(F_w) \times GL_h(F_w) \times \Q_p^\times) \lra
\Groth(\Frob_w^\Z) \]
which sends $[\alpha \otimes \beta \otimes \gamma]$ to
\[ \sum_\phi \vol(D_{F_w,n-h}^\times/F_w^\times)^{-1} \tr \alpha
(\varphi_{\Spp_{n-h}(\phi )}) (\dim \beta^{\Iw_{h,w}}) \Bigl[ \rec\bigl(\phi^{-1}
|\,\,\,|_w^{\frac{1-n}{2}} (\gamma^{\Z_p^\times} \circ \norm_{F_w/E_u})^{-1}\bigr)\Bigr], \]
where the sum is over characters $\phi$ of $F_w^\times/\CO_{F,w}^\times$.
(We just took the $\Iw_{h,w}$-invariant part of $\Red^{(h)}_1$, which
is defined on p.182 of \cite{ht}.)

Then, because
\[ Y_{U,S} = \bigcup_{T \supset S} Y_{U,T}^0 \]
for each $U=U^p\times \Iw(m)$, we have equalities
\begin{align*}
\bigl[H(Y_{\Iw(m),S},\CL_\xi)\bigr] &= \sum_{T\supset S}(-1)^{(n-\# S)-(n-\# T)} \bigl[H_c(Y_{\Iw(m),T}^0,\CL_\xi)\bigr]\\
&= \sum_{T\supset S}(-1)^{(n-\# S)-(n-\# T)} \bigl[H_c(I_{\Iw(m)}^{(n-\# T)},\CL_\xi)\bigr].
\end{align*}
As there are $\binom{n-\# S}{h}$ subsets $T$ with $\# T=n-h$ and $T\supset S$, we 
have proved the following lemma.

\begin{lem} \label{cohom}
If $l \neq p$, then for every $S\subset \{1,..., n\}$ we have an equality
\[\bigl[H(Y_{\Iw(m),S},\CL_\xi)\bigr] = \sum_{h=0}^{n-\# S}(-1)^{n-\# S-h} \binom{n-\# S}{h} \Red^{(h)} \bigl[H^i(X,\CL_\xi)\bigr]^{U_p^w(m)} \]
in the Grothendieck group of admissible $G(\A^{\infty,p}) \times \Frob_w^\Z$-modules over $\Q_l^{ac}$.
\end{lem}

The main innovation of our work is the following proposition.

\begin{prop} \label{middledeg} Suppose that $l \neq p$ and $\imath:\Q_l^{ac} \iso \C$. Suppose also that $\pi$ is an irreducible
admissible representation of $G(\A^\infty)$ such that
\begin{itemize}
\item $\pi_p^{\Iw(m)} \neq (0)$,
\item $R_{\xi,l}^j(\pi) \neq (0)$ for some $j$,
\item and $\BC_\imath(\pi)$ is cuspidal.
\end{itemize}
Then 
\[ H^j(Y_{\Iw(m),S},\CL_\xi)[\pi^p] =(0) \]
for $j \neq n - \# S$.
\end{prop}

\pfbegin
Let $\BC_\imath(\pi)=(\psi,\Pi)$. Then $H^j(X,\CL_\xi)^{U_p^w(m)}[\pi^p]=(0)$ if
$j \neq n-1$, while the space
$H^{n-1}(X,\CL_\xi)^{U_p^w(m)}[\pi^p]^{U^p}$ is $\imath^{-1}(\Pi_w
\times \psi_u)$-isotypic 
for any open compact subgroup $U^p \subset G(\A^{\infty,p})$. (See lemmas \ref{mo}
and \ref{coh}.) Moreover $\Pi_{w}$ has an Iwahori fixed vector and
$\psi_u$ is unramified, so that
\[( \dim \Pi_{w}^{\Iw_{n,w}}) \bigl[H^{n-1}(X,\CL_\xi)^{U^w_p(m)}[\pi^p]^{U^p}
\bigr]= (\dim H^{n-1}(X,\CL_\xi)^{\Iw(m)}[\pi^p]^{U^p}) \bigl[\imath^{-1}(\Pi_{w} \otimes
\psi_u)\bigr],\]
and 
\[n(\dim \Pi_{w}^{\Iw_{n,w}})\bigl[H_c(I_{\Iw(m)}^{(h)},\CL_\xi)[\pi^p]^{U^p}
\bigr] = (\dim H^{n-1}(X,\CL_\xi)^{\Iw(m)}[\pi^p]^{U^p})\,\imath^{-1} \Red^{(h)}\bigl[\Pi_{w}
\otimes \psi_u\bigr].\]
Combining this with lemma \ref{cohom}, we get
\begin{align*}
& n(\dim \Pi_{w}^{\Iw_{n,w}}) \bigl[H(Y_{\Iw(m),S},\CL_\xi)[\pi^p]^{U^p}\bigr]\\
&= (\dim H^{n-1}(X,\CL_\xi)^{\Iw(m)}[\pi^p]^{U^p}) \sum_{h=0}^{n-\# S} (-1)^{n-\# S -h} \binom{n-\# S}{h} \imath^{-1} \Red^{(h)}\bigl[\Pi_{w} \otimes \psi_u\bigr]. 
\end{align*}
As $\Pi_{w}$ is tempered, it is a full normalised induction of the form
\[ \nind_{P(F_w)}^{GL_n(F_w)} (\Spp_{s_1}(\pi_1) \otimes \cdots \otimes \Spp_{s_t}(\pi_t)), \]
where $\pi_i$ is an irreducible cuspidal representation of $GL_{g_i}(F_w)$ and
$P$ is a parabolic subgroup of $GL_n$ with Levi component $GL_{s_1g_1} \times
\cdots \times GL_{s_tg_t}$. As $\Pi_{w}$ has an Iwahori fixed vector, we must have 
$g_i=1$ and $\pi_i$ unramified for all $i$. Note that, for this type of
representation (full induced from square integrables $\Spp_{s_i}(\pi_i)$ with
$\pi_i$ an unramified character of $F_w^\times$), 
\begin{align*}
& \dim \bigl( \nind_{P(F_w)}^{GL_n(F_w)}(\Spp_{s_1}(\pi_1) \otimes \cdots \otimes \Spp_{s_t}(\pi_t)) \bigr)^{\Iw_{n,w}}\\
&= \# P(k(w)) \backslash GL_n(k(w)) / B_n(k(w)) = \frac{n!}{\prod_j s_j!}.
\end{align*}
We can compute $\Red^{(h)}\bigl[\Pi_{w} \otimes \psi_u\bigr]$ using lemma I.3.9
of \cite{ht} (but note the typo there --- ``positive integers $h_1,...,
h_t$" should read ``non-negative integers $h_1,..., h_t$"). Putting
$V_i=\rec\bigl(\pi_i^{-1}|\,\,\,|_w^{\frac{1-n}{2}}
(\psi_u \circ \norm_{F_w/E_u})^{-1}\bigr)$, we see that
\begin{align*}
\Red^{(h)}\bigl[\Pi_{w} \otimes \psi_u\bigr] &= \sum_{i} \dim
\bigl(\nind_{P'(F_w)}^{GL_h(F_w)} (\Spp_{s_i+h-n}(\pi_i|\,\,\,|^{n-h})
\otimes \bigotimes_{j \neq i} \Spp_{s_j}(\pi_j)) \bigr)^{\Iw_{h,w}}
\bigl[V_i\bigr]\\ &= \sum_{i} \frac{h!}{(s_i+h-n)! \prod_{j \neq i} s_j!}
\bigl[V_i\bigr]
\end{align*}
where the sum runs only over those $i$ for which $s_i\geq n-h$, and $P'
\subset GL_h$ is a parabolic subgroup. Thus
\begin{align*}
& n \frac{n!}{\prod_j s_j!} \bigl[H(Y_{\Iw(m),S},\CL_\xi)[\pi^p]^{U^p}\bigr]\\
&=  D \sum_{h=0}^{n-\# S}
(-1)^{n-\# S -h} \binom{n-\# S}{h} \sum_{i: \, s_i
\geq n-h} \frac{h!}{(s_i+h-n)! \prod_{j \neq i} s_j!} \bigl[V_i\bigr]\\
&= D \sum_{i=1}^t
\frac{(n-\#S)!}{(s_i-\#S)! \prod_{j \neq i} s_j!} \sum_{h=n-s_i}^{n-\# S}
(-1)^{n-\# S -h} \binom{s_i-\# S}{h+s_i-n} \bigl[V_i\bigr]\\
&= D \sum_{i\colon s_i=\# S}
\frac{(n-\# S)!}{\prod_{j \neq i} s_j!} \bigl[V_i\bigr],
\end{align*}
where $D=\dim H^{n-1}(X,\CL_\xi)^{\Iw(m)}[\pi^p]^{U^p}$, and so
\[n \binom{n}{\# S} \bigl[H(Y_{\Iw(m),S},\CL_\xi)[\pi^p]^{U^p}\bigr] =
(\dim H^{n-1}(X,\CL_\xi)^{\Iw(m)}[\pi^p]^{U^p}) \sum_{s_i=\# S}
\bigl[V_i\bigr].\]

As $\Pi_{w}$ is tempered, $\rec\bigl(\Pi_{w}^\vee \otimes (\psi_u^\vee \circ
\norm_{F_w/E_u})|\det|^{\frac{1-n}{2}}\bigr)$ is pure of weight $m_\xi-2t_\xi+(n-1)$.
Hence 
\[V_i=\rec\bigl(\pi_i^{-1}|\,\,\,|_w^{\frac{1-\# S}{2}}(\psi_u \circ
\norm_{F_w/E_u})^{-1} | \,\,\,|_w^{\frac{\# S-n}{2}}\bigr)\]
is strictly pure of weight $m_\xi-2t_\xi+(n-\# S)$. The Weil conjectures then tell us that
\[ H^j(Y_{\Iw(m),S},\CL_\xi)[\pi^p]^{U^p}=(0) \]
for $j \neq n-\# S$ and the proposition follows.
\pfend

\begin{cor}  Suppose that $l=p$, that $\tau:W_0 \into \Q_l^{ac}$ over $\Z_p=\Z_l$ and that $\imath:\Q_l^{ac} \iso \C$. Suppose also that
$\pi$ is an irreducible admissible representation of $G(\A^\infty)$ such that
\begin{itemize}
\item $\pi_p^{\Iw(m)} \neq (0)$,
\item $R_{\xi,p}^j(\pi) \neq (0)$ for some $j$,
\item and $\BC_\imath(\pi)$ is cuspidal.
\end{itemize}
Then 
\[ a_\xi (H^{j+m_\xi}(\CA^{m_\xi}_{\Iw(m),S}/W_0) \otimes_{W_0, \tau} \Q_l^{ac})[\pi^p] =(0) \]
for $j \neq n - \# S$.
\end{cor}

\pfbegin
Choose a prime $l' \neq l$ and an isomorphism $\imath':\Q_{l'}^{ac} \iso \C$. Set
$\xi'=(\imath')^{-1} \imath \xi$ and $\pi'=(\imath')^{-1} \imath \pi$. For any open compact subgroup $U^p \subset
G(\A^{\infty,p})$ we have
\begin{align*} 
& \dim_{\Q_l^{ac}} a_\xi (H^{j+m_\xi}(\CA^{m_\xi}_{\Iw(m),S}/W_0) \otimes_{W_0, \tau} \Q_l^{ac})[\pi^p]^{U^p} \\
= & \dim_{\Q_{l'}^{ac}} a_{\xi'}
H^{j+m_{\xi'}}(\CA^{m_{\xi'}}_{\Iw(m),S}, \Q_{l'}^{ac})[(\pi')^p]^{U^p} \\
= & \dim_{\Q_{l'}^{ac}} H^{j}(Y_{\Iw(m),S}, \CL_{\xi'})[(\pi')^p]^{U^p}. 
\end{align*}
(Use the main theorems of \cite{kme} and \cite{gm}.) The corollary now follows from the proposition.
\pfend

\begin{cor} Suppose that $\imath:\Q_l^{ac} \iso \C$. Suppose also that
$\pi$ is an irreducible admissible representation of $G(\A^\infty)$ such that
\begin{itemize}
\item $\pi_p^{\Iw(m)} \neq (0)$,
\item $R_{\xi,p}^j(\pi) \neq (0)$ for some $j$,
\item and $\BC_\imath(\pi)$ is cuspidal.
\end{itemize}
Then $\WD (R_{\xi,l}^{n-1}(\pi)|_{\Gal(F_w^{ac}/F_w)})$ is pure of weight $m_\xi-2t_\xi+n-1$.
\end{cor}

\pfbegin
The spectral sequence of proposition \ref{specseq}
\[ E_1^{i,j}(\Iw(m),\xi)[\pi^p] \Rightarrow \WD\bigl(H^{i+j}(X_{\Iw(m)}, \CL_\xi)|_{\Gal(F_w^{ac}/F_w)}\bigr)[\pi^p] \]
degenerates at $E_1$, as $E_1^{i,j}(\Iw(m),\xi)[\pi^p] = (0)$ unless $i+j=n-1$ by the proposition \ref{middledeg} and the previous corollary. Therefore the abutment is pure of desired weight, and the description of $H^{i+j}(X_{\Iw(m)}, \CL_\xi)$ in section \ref{s2} tells us that
\[\WD\bigl(H^{i+j}(X_{\Iw(m)}, \CL_\xi)|_{\Gal(F_w^{ac}/F_w)}\bigr)[\pi^p]=\bigoplus_{\pi' \sim \pi} (\pi')^{\Iw(m)} \otimes \WD(R_{\xi,l}^{n-1}(\pi')|_{\Gal(F_w^{ac}/F_w)})\]
(where the sum is over irreducible admissible $\pi'$ with $\pi'_x \cong \pi_x$ for all but finitely many finite places $x$ of $\Q$), hence the corollary.
\pfend

We now conclude the proof of theorem \ref{mt2} and hence of \ref{mt}. We return
to the notation of \ref{mt2}. Recall that $L$ is an
imaginary CM field and that $\Pi$ is a cuspidal automorphic representation
of $GL_n(\A_L)$ such that
\begin{itemize}
\item $\Pi \circ c \cong \Pi^\vee$;
\item $\Pi_\infty$ has the same infinitesimal character as some algebraic
representation over $\C$ of the restriction of scalars from $L$ to $\Q$ of
$GL_n$;
\item and for some finite place $x$ of $L$ the representation $\Pi_x$ is
square integrable.
\end{itemize}
Recall also that $v$ is a place of $L$ above a rational prime $p$, that
$l$ is a second rational prime and that $\imath:\Q_l^{ac} \iso \C$.
Recall finally that $R_l(\Pi)$ is the
$l$-adic representation associated to $\Pi$.

Choose a quadratic CM extension $L'/L$ in which $v$ and $x$ split. Choose
places $v'|v$ and $x'|x$ of $L'$ with $v'\ndiv x'(x')^c$. Also choose an
imaginary quadratic field $E$ and a totally real field $F^+$ such that
\begin{itemize}
\item $[F^+:\Q]$ is even;
\item $F=EF^+$ is soluble and Galois over $L'$;
\item $p$ splits as $p=uu^c$ in $E$;
\item if we denote by $\Pi_F$ the base change of $\Pi$ to $GL_n(\A_F)$, there is a place $w$ of $F$ above $u$ and $v'$ such that $\Pi_{F,w}$
has an Iwahori fixed vector;
\item $x$ lies above a rational prime which splits in $E$ and $x'$ splits in
$F$.
\end{itemize}
Note that the
component of $\Pi_F$ at a place above $x'$ is square integrable and hence
$\Pi_F$ is cuspidal.

Choose a place $z$ of $F$ above $x'$ and a division algebra $B$ with centre
$F$ as in section \ref{s2} and satisfying 
\begin{itemize}
\item $B_y$ is split for all places $y \neq z, z^c$ of $F$.
\end{itemize}
Also choose $\ast$, $\beta$ and $G$ as in section \ref{s2}. Then it follows
from theorem VI.2.9 and lemma VI.2.10 of \cite{ht} that we can find
\begin{itemize}
\item a character $\psi: \A_E^\times/E^\times \ra \C^\times$,
\item an irreducible algebraic representation $\xi$ of $G$ over $\Q_l^{ac}$,
\item and an irreducible admissible representation $\pi$ of $G(\A^\infty)$,
\end{itemize}
such that
\begin{itemize}
\item $R_{\xi,l}^{n-1}(\pi) \neq (0)$,
\item $\psi$ is unramified above $p$,
\item $\psi^c|_{E_\infty^\times}$ is the inverse of the restriction of
$\imath \xi$ to $E_\infty^\times \subset G(\R)$,
\item $\psi^c/\psi$ is the restriction of the central character of $\Pi_F$ to
$\A_E^\times$,
\item $\BC_\imath(\pi)=(\psi,\Pi_F)$,
\item and $\pi_w=\Pi_{F,w}$.
\end{itemize}

By the previous corollary $\WD((R_{\xi,l}^{n-1}(\pi)\otimes R_{l,\imath}(\psi)^{-1})|_{\Gal(F_w^{ac}/F_w)})$ is pure.
Moreover
\[ R_{\xi,l}^{n-1}(\pi)^\semis \otimes R_{l,\imath}(\psi)^{-1} \cong R_{l,\imath}(\Pi)
|_{\Gal(F^{ac}/F)}^a \]
for some $a \in \Z_{>0}$. Hence theorem \ref{mt2} follows from lemma \ref{gpure}.

\end{document}